\documentclass[10pt]{article}

\usepackage{amsmath,amssymb,stmaryrd}
\usepackage[all]{xy}

\begin{document}

\newtheorem{theorem}{Theorem}[section]
\newtheorem{lemma}[theorem]{Lemma}
\newtheorem{definition}[theorem]{Definition}
\newtheorem{example}[theorem]{Example}
\newtheorem{corollary}[theorem]{Corollary}

\title{Strong phase-space semiclassical asymptotics}

\author{
Agissilaos G. ATHANASSOULIS%
\footnote{CMLS, \'Ecole Polytechnique, 91128 Palaiseau,
athanas@math.polytechnique.fr},
Thierry PAUL 
\footnote{ 
CNRS, CMLS, \'Ecole Polytechnique, 91128 Palaiseau, 
thierry.paul@math.polytechnique.fr}
}
\date{}

\maketitle

%

\begin{abstract}
Wigner and Husimi transforms have long been used for the phase-space reformulation of Schr\"odinger-type equations, and the study of the corresponding semiclassical limits. Most of the existing results provide approximations in  appropriate weak topologies.
In this work we are concerned with semiclassical limits in the strong topology, i.e. approximation of Wigner functions by solutions of the Liouville equation in $L^2$ and Sobolev norms. The results obtained improve the state of the art, and highlight the role of potential regularity, especially through the regularity of the Wigner equation. It must be mentioned that the strong convergence can be shown up to $O(log \frac{1}\varepsilon)$ time-scales, which is well known to be, in general, the limit of validity of semiclassical asymptotics.
\end{abstract}

\setcounter{tocdepth}{2}
\tableofcontents


\section{Introduction and main results}

Consider a wavefunction $u^\varepsilon(x,t)$ satisfying the Schr\"{o}dinger equation,
\begin{equation}
\label{eq1}
\begin{array}{c}
i\varepsilon \frac{\partial }{\partial t}u^\varepsilon = \left( -\frac{ \varepsilon^2}{2}\Delta +V\left( x \right) \right) u^\varepsilon , \\
u^\varepsilon(t=0)=u^\varepsilon_0(x).
\end{array}
\end{equation}

The Wigner transform (WT) $W^\varepsilon(x,k,t)$ of the wavefunction $u$ is defined in the standard semiclassical scaling,
\begin{equation}
 W(x,k,t)=\int\limits_{y \in \mathbb{R}^n} { e^{-2\pi i yk}u^\varepsilon (x+\varepsilon \frac{y}2,t) \bar{u}^\varepsilon (x-\varepsilon \frac{y}2,t) dy }.
\end{equation}
In that case $W$ corresponds to a {\em pure state}.
A Wigner function $W$ need not be the Wigner transform of a wavefunction $u$; when working with {\em mixed states} the Wigner function can be
\begin{equation}
\label{eqdefrelrowig}
 W^\varepsilon(x,k,t)=\int\limits_{y \in \mathbb{R}^n} { e^{-2\pi i yk} \rho^\varepsilon(x+\varepsilon \frac{y}2,x-\varepsilon \frac{y}2,t) dy }
\end{equation}
for the {\em density matrix} $\rho^\varepsilon(x,y)$, which satisfies the von Neumann-Heisenberg equation
\begin{equation}
\begin{array}{c}
\label{eqtralala}
i \varepsilon \rho^\varepsilon_t + \frac{\varepsilon^2}{2}(\Delta_x-\Delta_y)\rho^\varepsilon - (V(x)-V(y))\rho^\varepsilon =0, \\ { } \\
\rho^\varepsilon(t=0)=\rho_0^\varepsilon.
\end{array}
\end{equation}
($\rho_0^\varepsilon$ corresponds to a pure state if there is a $u^\varepsilon_0(x)\in L^2$ such that $\rho_0^\varepsilon(x,y)=u_0^\varepsilon(x)\overline{u}^\varepsilon_0(y)$).

Now denote by $\Phi$ the $\varepsilon$-dependent smoothing operator
\begin{equation}
\label{eqDefPhi}
\begin{array}{c}
 \Phi:f(x,k) \mapsto \left({ \frac{2}{\varepsilon \sigma_x\sigma_k} }\right)^n \int{e^{-\frac{2\pi}{\varepsilon} \left[{ \frac{|x-x'|^2}{\sigma_x^2} + \frac{|k-k'|^2}{\sigma_k^2} }\right] }  f(x',k') dx'dk'}.
\end{array}
\end{equation}
$\sigma_x$ and $\sigma_k$ are $\varepsilon$-independent parameters that can be used for fine tuning. 
As a matter of notation we will use $\widetilde{f}=\Phi f$.

Whether we deal with pure or mixed states, the smoothed Wigner transform (SWT) is defined as
\begin{equation}
 \widetilde{W}=\Phi W.
\end{equation}

The WT and SWT satisfy equations associated with (\ref{eq1}), (\ref{eqtralala}). As was shown in \cite{AMP}, the equation for the SWT $\widetilde{W}^\varepsilon(x,k,t)$ is
\begin{equation}
\label{SWTeq}
\begin{array}{l}
\partial_t \widetilde{W}^\varepsilon + \left({ 2\pi k \cdot \partial_x +\frac{\varepsilon\sigma_x^2}{2}\partial_x\cdot\partial_k }\right) \widetilde{W}^\varepsilon+ \\ 

\,\,\,\,\,\,\,\,\,\,\,\,\,\,\,\,\,
+\frac{2}{\varepsilon} Re \left[{ 
i\int{e^{2\pi i Sx-\frac{\varepsilon\pi}{2}\sigma_x^2S^2}\hat{V}(S)\widetilde{W}^\varepsilon(x+\frac{i\varepsilon\sigma_x^2S}{2},k-\frac{\varepsilon S}{2})dS}  }\right]=0, \\ { } \\
\,\,\,\,\,\,\,\,\,\,\,\,\,\,\,\,\,\,\,\,\,\,\,\,\,\,\,\,\,\,\,\,\,\,\,\,\,\,\,\,\,\,\,\,\,\,\,\,\,\,\,\,\,\,\,\,\,\,\,\,\,\,\,\,\,\,
\widetilde{W}(t=0)=\widetilde{W}^\varepsilon_0.
\end{array}
\end{equation}
In accordance with the meaning of $\sigma_x,\sigma_k$, the well-known equation for $W^\varepsilon$ can be recovered by setting $\sigma_x=\sigma_k=0$ in (\ref{SWTeq}).

It is well known that given a family of problems as in (\ref{eq1}), studying the limit $\varepsilon \shortrightarrow 0$  presents difficulties, both in terms of analysis (e.g. the appearence of caustics; also worth mentioning that typically $u^\varepsilon$ doesn't have a meaningful limit in $\varepsilon$) as well as computation (the complexity of solving numerically problem (\ref{eq1}) blows up as $\varepsilon \shortrightarrow 0$). That's the reason why several asymptotic techniques have been developed to treat semiclassical (also called {\em high-frequency}, or {\em short-wavelength} in some contexts) problems -- one approach being based on Wigner measures (WMs).

WMs are a very powerful and elegant tool, which under minimal assumptions in many cases allows a simple and elegant description of the semiclassical limit. On a mathematical level, the idea is that (in an appropriate sense) $\mathop{lim}\limits_{\varepsilon \shortrightarrow 0} W^\varepsilon$ is a natural object (even when $\mathop{lim}\limits_{\varepsilon \shortrightarrow 0} u^\varepsilon$ is not), and evolves in time under simple asymptotic dynamics.

A very simple way to state the core result is that, if the potential $V$ is smooth enough then (up to the extraction of a subsequence)
\[
W^\varepsilon(t) \rightharpoonup W^0(t)
\] 
in an appropriate weak-$*$ sense, where
\[
\begin{array}{c}
\partial_t W^0 +2\pi k \cdot \partial_k W^0 - \frac{1}{2\pi} \partial_x V \cdot \partial_k W^0=0,\\ { } \\
W^0(t=0)=\mathop{lim}\limits_{\varepsilon \shortrightarrow 0} W^{\varepsilon}_0,
\end{array}
\]
see e.g. \cite{LP}. The Liouville equation can then be solved with the method of characteristics, i.e. is effectively reduced to an ODE problem.

WMs have been also used to derive asymptotic models in problems featuring systems of equations, nonlinearities, stochastic or periodic coefficients, inverse problems etc; see for example \cite{LP,ee,ey1,GMMP,Mark,MM,MMP,LM,Pap,RPK} and the references therein.
However, there are intrinsic limitations to the WM approach, that have only recently started to attract systematic attention -- see e.g. \cite{cfk}, section \ref{ss003}.

\vskip 0.7cm
The motivation of the present work is threefold, and can be summarized as follows: 

\vskip 0.2cm
{\em On a technical-mathematical level}, we show results that improve upon the state of the art on strong-topology 
semiclassical limits. See section \ref{ss001} for more details. It must be noted that the application (and extension) of the machinery  developed in \cite{AMP} is crucial to the proofs here. Moreover, all the error estimates are constructive, most of the constants are computed explicitly -- in fact all of them are computable, including in particular the time dependence.

\vskip 0.2cm
{\em On a numerical level}, our results justify the use of the SWT for fast coarse-scale simulations; see e.g. \cite{a1,a2}, and section \ref{ss002}.

\vskip 0.2cm
Finally, {\em on a more qualitative level,} the ideas developed here can be used as an intermediate step to work out problems in which the semiclassical limit is not effectively known; see section \ref{ss003}.

\subsection{Statement of the assumptions}
Let $r \in \mathbb{N} \cup \{ 0 \}$ be fixed. The main assumption will be

\noindent  {\bf Assumption A1(r)} {\em There is a $M_0>0$ such that
\begin{equation}
\label{eqmklop00}
\begin{array}{c}
 \int{|\hat{V}(S)|(1+|S|^{r+2})dS}=M_0.
\end{array}
\end{equation}

Now denote the auxiliary potentials $V_1(x)$, $\widetilde{V}_1(x)$ to be
\begin{equation}
\label{eqV1}
\begin{array}{c}
V_1(x) = \left({ \frac{\pi}{\eta} }\right)^{\frac{n}{2}} \int{e^{-\frac{\pi^2}{\eta}|x-x'|^2} V(x')dx' }, \\ { } \\

\widetilde{V}_1(x) = \left({ \frac{2}{\varepsilon\sigma_x^2} }\right)^{\frac{n}{2}} \int{e^{-\frac{2\pi}{\varepsilon\sigma_x^2}|x-x'|^2} V_1(x')dx' }.
\end{array}
\end{equation}
The parameter $\eta$ is set to be
\[
\eta =  \left({  \frac{\pi}2 \cdot max\left({   (4n-1)\sigma_x^2 \,\, , \,\,   \frac{4n}{\sigma_k^2}-\sigma_x^2   }\right) +1 }\right)  \varepsilon
\]
 (see equation (\ref{SetEta}) and above for the rationale behind this scaling). 

\vskip 0.2cm
There are $C,D,R>0$, $\theta \in (r+1,r+2)$ such that

\begin{equation}
\label{eqq01}
\begin{array}{c}
\forall m\in \mathbb{N}, m < \theta +1, \,\,\,:\,\,\, |A|\leqslant m \,\,\,\Rightarrow \,\,\,
||\partial_x^{A} \widetilde{V}_1(x)||_{L^\infty} \leqslant  C,
\end{array}
\end{equation}
\begin{equation}
\label{eqq02}
\begin{array}{c}
\forall m\in \mathbb{N}, m > \theta +1, \,\,\,:\,\,\, |A|\leqslant m \,\,\,\Rightarrow \,\,\,
||\partial_x^{A} \widetilde{V}_1(x)||_{L^\infty} \leqslant  \\ { } \\
\leqslant   (2\pi)^m M_0 +D \frac{(2\pi)^m}{2} (\eta+\varepsilon\frac{\pi}{2}\sigma_x^2)^{-\frac{m-1-\theta}{2}} \Gamma\left({ \frac{m-1-\theta}{2}}\right).
\end{array}
\end{equation}

In addition, $\forall m =0,1,..,r$
\begin{equation}
\label{eqq03}
|| (1-e^{-(\eta+\varepsilon\frac{\pi}{2}\sigma_x^2) S^2}) |S|^{m+1}\hat{V}(S)  ||_{L^1(\mathbb{R}^{n})} \leqslant C(\varepsilon^\frac{\theta-m}{2}+\varepsilon)
\end{equation}
}

\vskip 0.5cm
There is a somewhat stronger version of this assumption which is much simpler:

\noindent  {\bf Assumption A1'(r)} {\em There is a $M_0>0$ such that
\begin{equation}
\label{eqmklop007}
\begin{array}{c}
 \int{|\hat{V}(S)|dS}=M_0,
\end{array}
\end{equation}
and
$\exists R,D>0$, $\theta \in (r+1,r+2)$ such that
\begin{equation}
\label{eqmklop}
\begin{array}{c} 
\forall |S|>R \,\,:\,\, |\widehat{V}(S)|\leqslant D |S|^{-(n+1+\theta)}.
\end{array}
\end{equation}
}

\vskip 0.3cm
We mentioned that $A1'(r) \,\, \Rightarrow \,\, A1(r)$. Indeed,

\begin{lemma} \label{lemmqq} Equations (\ref{eqmklop007}) and (\ref{eqmklop}) imply
(\ref{eqmklop00}), (\ref{eqq01}), (\ref{eqq02}), (\ref{eqq03}).
\end{lemma}

\vskip 0.15cm
\noindent {\bf Proof:} The proof is given at the beginning of section \ref{ecccqq}.

\vskip 0.25cm
In some cases we will use 

\noindent {\bf Assumption A2} {\em  A1(r=0) holds, and moreover $||W^\varepsilon_0||_{H^1}= o(\varepsilon^{-\frac{1}{2}} ||W_0^\varepsilon||_{L^2})$.
}

\vskip 0.2cm
\noindent {\bf Remarks:}
\vskip 0.25cm
\noindent {\bf $\bullet$} The reason we use the more complicated (but weaker) conditions (\ref{eqq01}), (\ref{eqq02}), (\ref{eqq03}) instead of (\ref{eqmklop}) is that there are relevant examples not covered by it;  a simple example is given by 
a lower-dimensional potential: consider $V_*(x_1,...,x_n)=V_*(x_1,...,x_d)$ with $d<n$; if $V_*$ satisfies
\[
\forall |S|>R \,\,:\,\, |\widehat{V}_*(S_1,...,S_d)|\leqslant D \left({\sqrt{\sum\limits_{j=1}^d S_j^2} }\right)^{-(d+1+\theta)}
\]
it is easy to check that A1 follows -- although in general (\ref{eqmklop}) does not.

\noindent {\bf $\bullet$}  Following assumption A1(r=0), $V(x) \in L^\infty$. This suffices to imply that
the Schr\"odinger operators $-\frac{1}2\Delta +V(x)$, $-\frac{1}2\Delta +V_1(x)$, $-\frac{1}2\Delta +\widetilde{V}_1(x)$  are essentially self-adjoint on $L^2(\mathbb{R}^n)$. This is something that we will use, e.g. in applying theorem \ref{thrmWignL2} (and the reason why we introduce $V_1$ here).
 
We will not comment on this explicitly any more.

\vskip 0.25cm
\noindent {\bf $\bullet$} It is clear that if A1(r) is true, then $\forall  r'\in [0,r] \cap \mathbb{Z} \,\, A1(r')$ is also true.

\subsection{Formulation of the main results} \label{sbsc007}
\begin{theorem}[$L^2$ semiclassical asymptotics] \label{thrm2}

Consider a potential $V(x)$ and an initial Wigner function $W_0^\varepsilon$, such that (A2) is satisfied. 

 Denote by $W^\varepsilon(x,k,t)$ the solution of the corresponding Wigner equation with initial data $W^\varepsilon_0$ as in (A2), $\widetilde{W}^\varepsilon(t)=\Phi W^\varepsilon(t)$ the corresponding SWT, $\rho(x,k,t)$ the solution of the Liouville equation
\begin{equation}
\label{2LiouvEqa}
\begin{array}{c}
\partial_t \rho +2\pi k \cdot \partial_k \rho - \frac{1}{2\pi} \partial_x V \cdot \partial_k \rho=0,\\ { } \\
\rho(t=0)=W^{\varepsilon}_0,
\end{array}
\end{equation}
and by $\rho_1^\varepsilon(x,k,t)$ the solution of
\begin{equation}
\label{SmLiouvEqa}
\begin{array}{c}
\partial_t \rho_1^\varepsilon +2\pi k \cdot \partial_k \rho_1^\varepsilon - \frac{1}{2\pi} \partial_x \widetilde{V}_1 \cdot \partial_k \rho_1^\varepsilon=0,\\ { } \\
\rho_1^\varepsilon(t=0)=\widetilde{W}^\varepsilon_0.
\end{array}
\end{equation}

Then there is an $O(1)$ explicit constant $C$, depending only on $n,\sigma_x,\sigma_k,M_0$ such that, for all $t \in [0,T]$, $m \in \mathbb{N}$,
\begin{equation}
\label{equjmn1}
\begin{array}{c}
|| \rho_1^\varepsilon(t)-\widetilde{W}^\varepsilon(t) ||_{L^2} \leqslant \\ { } \\ 

\leqslant \sqrt{\varepsilon}  CT e^{Tn \pi \, max \{ 1, ||\hat{V}(S)|S|^{2}||_{L^1} \} } ||W^\varepsilon_0||_{H^1},
\end{array}
\end{equation}

\begin{equation}
\label{eqaqjmn1}
\begin{array}{c}
\frac{
|| \rho_1^\varepsilon(t)-\widetilde{W}^\varepsilon(t) ||_{L^2} 
}{
||\widetilde{W}_0^\varepsilon ||_{L^2} 
} 

\leqslant C \sqrt{\varepsilon} T e^{Tn \pi \, max \{ 1, ||\hat{V}(S)|S|^{2}||_{L^1} \} } \frac{||W^\varepsilon_0||_{H^1}}{||W^\varepsilon_0||_{L^2}}=\\ { } \\

= o\left( { T e^{Tn \pi \, max \{ 1, ||\hat{V}(S)|S|^{2}||_{L^1} \} } }\right) ,
\end{array}
\end{equation}

\begin{equation}
\label{equjmn2}
\begin{array}{c}
||\rho(t)-W^\varepsilon(t) ||_{L^2} \leqslant \\ { } \\

\leqslant C  \left( { \sqrt{\varepsilon}  (T+1)||W^\varepsilon_0||_{H^1}+T\varepsilon ||W^\varepsilon_0||_{L^2} }\right) e^{Tn \pi \, max \{ 1, ||\hat{V}(S)|S|^{2}||_{L^1} \} },
\end{array}
\end{equation}

\begin{equation}
\label{eqaqjmn2}
\begin{array}{c}
\frac{||\rho(t)-W^\varepsilon(t) ||_{L^2}}{||W^\varepsilon_0 ||_{L^2}} 

\leqslant C \sqrt{\varepsilon} (1+T)e^{Tn \pi \, max \{ 1, ||\hat{V}(S)|S|^{2}||_{L^1} \} } \frac{||W^\varepsilon_0||_{H^1}}{||W^\varepsilon_0||_{L^2}}= \\ { } \\

=o\left( {  (1+T)e^{Tn \pi \, max \{ 1, ||\hat{V}(S)|S|^{2}||_{L^1} \} }  }\right).
\end{array}
\end{equation}
\end{theorem}

\noindent {\bf Remarks:} 

\vskip 0.25cm
\noindent {\bf $\bullet$} Observe that the result can be extended beyond time-scales $T=O(1)$, up to $T=O(-log(\varepsilon))$. As is well known, semiclassical asymptotics in general cannot be extended for longer time-scales \cite{Rob,Paul}. It is worth noting that we can go to such time scales with A1(r=0), which amounts to 2 derivatives on the potential, as opposed to analytic potentials (in \cite{BGP}) or $C^{\infty}$ potentials (in \cite{Rob}).

\vskip 0.25cm
\noindent {\bf $\bullet$} Similar approximations hold e.g. between $W^\varepsilon(t) \approx \rho^\varepsilon_1(t)$, etc; see the proof. 

\vskip 0.25cm
\noindent {\bf $\bullet$} One can see in the proof that all constants are (or can easily be) computed explicitly. Although of course they will be pessimistic, given a problem (i.e. potential and initial data) they provide one with an a priori estimation of all the errors involved. This is true for the following theorems as well.

\vskip 0.25cm
\noindent {\bf $\bullet$} Pure states, i.e. $W^\varepsilon_0=W^\varepsilon[u^\varepsilon]$, typically give $\frac{||W^\varepsilon_0||_{H^1}}{||W^\varepsilon_0||_{L^2}} =O(\varepsilon^{-\frac{n}{2}})$ -- see lemma \ref{lemalan}. Therefore this result does not yield a significant estimate for pure states (i.e. the relative error bound obtained is not $o(1)$). See corollary \ref{cor100} and the remark thereafter for pure states.

\vskip 0.4cm
\begin{theorem}[Negative-index Sobolev spaces]
\label{thrm23}
Consider some $r \in \mathbb{N}$. Assume that A1(r) is satisfied, and $W_0^\varepsilon \in H^{1-r}$.

Then,  for all  $t\in [0,T]$ 
\begin{equation}
\begin{array}{c}
\label{eqaqjmn1NEG}
|| \rho^\varepsilon(t)-{W}^\varepsilon(t) ||_{H^{-r}}  \leqslant  C T ||W_0^\varepsilon||_{H^{-1-r}} e^{TD(r+1,n)} (\sqrt{\varepsilon} + \varepsilon^{\frac{\theta}{2}}) + \\ { } \\ + C (e^{TD(r,n)}+e^{T[D(r,n)+D(r-1,n)]}) \sqrt{\varepsilon} ||W_0^\varepsilon||_{H^{1-r}}+ \\ { } \\
+ (\varepsilon+\varepsilon^{\frac{\theta}{2}})
C T e^{T[D(r+1,n)+D(r,n)]} || W_0^\varepsilon||_{H^{-1-r}} .
\end{array}
\end{equation}

As long as there is a constant $C_1$, independent of $\varepsilon$, so that $||W^\varepsilon(t)||_{H^{-r}} \geqslant C_1$, equation (\ref{eqaqjmn1NEG}) yields a relative error estimate in a straightforward manner.
\end{theorem}

\vskip 0.25cm
\noindent {\bf Remark: } As long as the assumtpions of \cite{LP} are satisfied, it can be used to bound from below
$||W^\varepsilon(t)||_{H^{-r}}$. Take $\phi \in \mathcal{S}$, $||\phi||_{H^r} \leqslant 1$; \cite{LP} assures that $\langle W^\varepsilon,\phi \rangle \shortrightarrow \langle W^0,\phi \rangle$. So a sufficient condition is that there is a $\phi_0(t)$ so that, for $t \in [0,T]$, $|\langle W^0(t),\phi_0(t) \rangle | > 2C_1 >0$. 

If that's true, then it follows automatically that 
\[
\exists \varepsilon_0 : \forall \varepsilon \in (0,\varepsilon_0): | \langle W^0(t) - W^\varepsilon(t),\phi_0(t) \rangle |<C_1.
\]

Of course this is mostly a technical point; we are not aware of any meaningful example where $||W^\varepsilon_0||_{H^{-r}}=O(1)$, and later $||W^\varepsilon(t)||_{H^{-r}}=o(1)$.

\vskip 0.4cm
\begin{theorem}[Positive-index Sobolev spaces]
\label{thrm24}
Consider some $r \in \mathbb{N}$. Assume that A1(r) is satisfied, and $W_0^\varepsilon \in H^{1+r}$.

Then,  for all  $t\in [0,T]$ 
\begin{equation}
\begin{array}{c}
\label{eqaqjmn1POS}
|| \rho^\varepsilon(t)-{W}^\varepsilon(t) ||_{H^{r}}  \leqslant  C || {W}_0^\varepsilon ||_{H^{1+r}} 
\left[{ Te^{TD(r+1,n)}(\sqrt{\varepsilon}+\varepsilon^{\frac{\theta}2}) +}\right. \\ { } \\

\left.{+ Te^{T[D(r+1,n)+D(r,n)]}(\varepsilon+\varepsilon^{\frac{\theta-r}2}) + e^{TD(r,n)\sqrt{\varepsilon}}

}\right]
\end{array}
\end{equation}

Moreover, if in addition
\begin{equation}
\label{eqaqjmn1POSrel}
\frac{|| {W}_0^\varepsilon ||_{H^{1+r}} }{|| {W}^\varepsilon(t) ||_{H^{r+1}}} =o(\varepsilon^{-\frac{1}{2}}),
\end{equation}
equation (\ref{eqaqjmn1POS}) yields directly a relative error estimate as well.
\end{theorem}

\noindent {\bf Remark: } A sufficient condition for equation (\ref{eqaqjmn1POSrel}) is given by 
$ \frac{ ||W^\varepsilon_0||_{H^{r+1}}  }
{||W^\varepsilon_0||_{L^2}} =o(\varepsilon^{-\frac{1}{2}})$. Indeed, to see that, apply
\[
|| {W}^\varepsilon(t) ||_{H^{r+1}} \geqslant || {W}^\varepsilon(t) ||_{L^2} \,\,\, \Rightarrow \,\,\,
\frac{1}{|| {W}^\varepsilon(t) ||_{H^{r+1}}} \leqslant \frac{1}{|| {W}^\varepsilon(t) ||_{L^2}}=\frac{1}{|| {W}^\varepsilon_0 ||_{L^2}} 
\]
to estimate $\frac{1}{|| {W}^\varepsilon(t) ||_{H^{r+1}}}$ in the lhs of equation (\ref{eqaqjmn1POSrel}).

\vskip 0.8cm
It is interesting to look at some corollaries:
\begin{corollary}[$L^2$ limit]\label{cor1} Assume that  (A2) holds, and  $\exists W^0_0 \in L^2(\mathbb{R}^{2n})$ such that
\[
\mathop{lim}\limits_{\varepsilon \shortrightarrow 0} ||W^0_0- W^\varepsilon_0 ||_{L^2}=0.
\]
Denote $\rho^0(x,k,t)$ the solution of
\begin{equation}
\label{2LiouvEqa2}
\begin{array}{c}
\partial_t \rho^0 +2\pi k \cdot \partial_k \rho^0 - \frac{1}{2\pi} \partial_x V \cdot \partial_k \rho^0=0,\\ { } \\
\rho^0(t=0)=W^{0}_0.
\end{array}
\end{equation}

Then, $\forall t \in [0,T]$,
\[
\mathop{lim}\limits_{\varepsilon \shortrightarrow 0} ||\rho^0(t)- W^\varepsilon(t) ||_{L^2}=0.
\]
\end{corollary}

\vskip 0.4cm
\begin{corollary}[Strong limit with concentration]\label{cor100} Assume that  $\exists s \in \mathbb{N}, W^0_0 \in H^{1-s}$ such that
\[
\mathop{lim}\limits_{\varepsilon \shortrightarrow 0} ||W^0_0 - W^\varepsilon_0||_{H^{-s}} =0.
\]
Assume moreover that assumption A1(s) is satisfied, and
recall the definition of  $\rho^0(x,k,t)$,  in equation (\ref{2LiouvEqa2}).

Then, $\forall t \in [0,T]$,
\[
\mathop{lim}\limits_{\varepsilon \shortrightarrow 0} ||\rho^0(t) - W^\varepsilon(t)||_{H^{-s}} =0.
\]
\end{corollary}

\noindent {\bf Remark: Pure states.}  It is important to note that this result can also be applied to {\em pure states}: for example if the initial data is a coherent state, then (provided the potential is smooth enough) the assumptions of corollary \ref{cor100} are satisfied for $s=\lceil{ \frac{n}{2} }\rceil+2$; see lemma \ref{LEMATARA}.

%


\vskip 0.4cm
However, strong aproximation doesn't have to mean strong limit: 
\begin{corollary}[Strong approximation with concentration]\label{cor2} Assume (A2) holds, and $W^\varepsilon_0 \rightharpoonup W^0_0 \notin L^2$ in weak-$*$ sense (as in \cite{LP}). 

Then
\[
\frac{ || {W}^\varepsilon(t) - \rho^\varepsilon(t) ||_{L^2} }{ ||{W}^\varepsilon_0||_{L^2} }=
o(1).
\]
\end{corollary}
\noindent {\bf Remark: } 
While neither $\rho^\varepsilon(t)$, nor ${W}^\varepsilon(t)$ have a limit in $L^2$ as $\varepsilon \rightarrow 0$, they are asymptotically close to each other in $L^2$ $\forall \varepsilon >0$. The point here is  that there are situations where this might apply while corollary \ref{cor100} not -- for example when the convergence to the Wigner measure isn't known in norm, or when the potential isn't smooth enough.

Observe finally that an $L^2$ approximation is much stronger than and $H^{-s}$ limit. For example, the $L^2$ approximation implies any oscillations that might develop during the quantum evolution (i.e. absent from the initial data) are necessarily small in $L^2$ sense -- while of course the $H^{-(\lceil \frac{n}{2} \rceil+2)}$ limit tells us only that eventually any oscillations will cancel out in weak sense.

\vskip 0.7cm
\noindent  Finally there are a couple of somewhat technical points worth mentioning:

\vskip 0.25cm
\noindent {\bf $\bullet$} As long as $\sigma_x \sigma_k \geqslant 1$, the SWT ``preserves positivity'', i.e. if $W^\varepsilon(x,k)$ is the Weyl symbol of a positive operator, then $\widetilde{W}^\varepsilon(x,k) \geqslant 0$. 
(This is well-known, and is a property that  the WT itself does not have; see lemma \ref{lemhuscrit}). Our asymptotic analysis here does not assume anything on the sign of $\widetilde{W}^\varepsilon(x,k)$.  Still, it should be noted that a strong approximation which preserves sign might be of particular interest in more difficult (e.g. nonlinear) problems.

\vskip 0.25cm
\noindent {\bf $\bullet$} It is also important to note that the basic limitations come from the breakdown of Sobolev regularity in $\rho_1^\varepsilon$, $W^\varepsilon$. Therefore, to understand where there is an actual breakdown of the approximation, as opposed to a mere technical limitation of our tools, it will be important to get sharper assumptions on the validitiy -- and failure -- of theorems such as \ref{thrm1corr}. The weaker assumptions e.g. of theorem \ref{thrmnext} is one possible way forward.

\subsubsection{State of the art}
\label{ss001}

One reason for working out strong topology semiclassical approximations for the Wigner function, is that
 there is an extremely limited literature on the subject. Indeed, most results focus on weak-$*$ limits.
 
 From an analytical point of view, a norm could quantify a rate of convergence -- in contrast to the weak-$*$ convergence, which depends on the choice of test function. This could also be extremely important  with respect to computational applications (see next section), where it is important to quantify a priori (even if on the pessimistic side) the errors involved. From a more qualitative point of view, it could be the way  to  add $\varepsilon$-dependent corrections to the limit, thus controlling in a more flexible way the information loss.

  As long as we are working with norms, $L^2$ norms are especially well suited -- in particular in view of theorem \ref{thrmWignL2}. The $H^m$ generalizations  also, in fact,
make use of theorem \ref{thrmWignL2}. ($L^1$ approximations also make sense, in view of the classical interpretation. We don't work on that here, but it may be a direction worth pursuing in the future).

 To the best of our knowledge  \cite{Pul} comprises the state of the art for strong semiclassical approximations. Other existing works (e.g. \cite{BGP}) often use regularity assumptions stronger than \cite{Pul}, and, to our knowledge, do not treat $H^m$ norms. It is therefore natural that our results here should be compared to those of \cite{Pul}.

In our notation, Theorem 5.1 of \cite{Pul} states (or rather implies, since it also treats higher order approximations):

\noindent {\bf Statement}{\em \,\,
If
\begin{equation}
\label{eqtgfolk}
||W_0^\varepsilon||_{H^2} < \infty, \,\,\,\,\, \int{|\widehat{V}(S)|(1+|S|^2)dS}<\infty,
\end{equation}
then
\begin{equation}
||\rho(t)-W^\varepsilon(t)||_{L^2} \leqslant O\left({\,\, \varepsilon^2  ||W^\varepsilon_0||_{H^2} e^{Ct(1+\int{|\widehat{V}(S)|\,|S|^2dS})} \,\,}\right).
\end{equation}
}

\vskip 0.25cm
However, there is a misprint in \cite{Pul}; the statement of the theorem should read
\begin{theorem}
If
\begin{equation}
\label{eqtgfolkoo}
||W_0^\varepsilon||_{H^2} < \infty, \,\,\,\,\, \int{|\widehat{V}(S)|(1+|S|^3)dS}<\infty,
\end{equation}
then
\begin{equation}
||\rho(t)-W^\varepsilon(t)||_{L^2} \leqslant O\left({\,\, \varepsilon^2  ||W^\varepsilon_0||_{H^2} e^{Ct(1+\int{|\widehat{V}(S)|\,|S|^3dS})} \,\,}\right).
\end{equation}
\end{theorem}

(The problem in Theorem 5.1 originally comes from from a typo in Theorem 4.1. See section \ref{qrr}, theorem \ref{thrm1corr} for a detailed treatment of the point in question).

\vskip 0.25cm
In that light, the use of the SWT in Theorem \ref{thrm2},  allows one to gain one derivative in each of the initial data and the potential. (For example $\widehat{V}(S)=\frac{1}{(1+|S|)^{n+\frac{5}{2}}}$ satisfies (A2), while $\int{|\widehat{V}(S)| \,\, |S|^3dS}=\infty$).

The price one has to pay for these weaker regularity assumptions, is a larger size of the error: for example, if $\theta >2$ and $||W_0^\varepsilon||_{H^2} =O(1)$, theorem \ref{thrm2} gives $||\rho(t)-W^\varepsilon(t)||_{L^2} \leqslant O(\sqrt{\varepsilon})$, as opposed to an $O(\varepsilon^2)$ error from Theorem 5.1 of \cite{Pul}.

It must also be noted that neither result applies to pure states. \cite{Pul} contains a result in $H^{-m}$ developed with pure states in mind, Theorem 5.2; our approach here also yields an $H^{-m}$ approximation result applicable to pure states, corollary \ref{cor100} which follows from theorem \ref{thrm23}. (Note that theorem \ref{thrm23} also improves on the regularity assumptions of Theorem 5.2 in \cite{Pul}; we don't do the exhaustive comparison here).

\subsubsection{Numerical applications}
\label{ss002}

Part of the motivation for this work has been its potential for numerical applications. Indeed, a numerical scheme using the SWT for the fast recovery of coarse-scale observables was formulated in \cite{a1,a2}, and proof-of-concept numerical experiments showed that it is a competitive option in the semiclassical regime. This work serves as a rigorous justification for the method. 

\vskip 0.25cm
In this regard it is worth highlighting one point in particular:
it is well known that  keeping track of the WM allows one to approximate quadratic observables of interest, with famous example
\[
|u^\varepsilon(x,t)|^2 \approx  \int{W^0(x,k,t)dk}.
\]
In addition, WMs don't break down -- or even require any special treatment -- when caustic appear. However, although a WM stays well defined on phase-space, {\em its marginals need not exist, especially when a caustic appears}. That is, while 
\[
\forall \varepsilon >0 \,\,\,\,  |u^\varepsilon(x,t)|^2 = \int{W^\varepsilon[u^\varepsilon(t)](x,k,t)dk},
\] 
$\int{W^0(x,k,t)dk}$ may not be well defined, even if $\mathop{lim}\limits_{\varepsilon \shortrightarrow 0} |u^\varepsilon(t)|^2$ is\footnote{$\mathop{lim}\limits_{\varepsilon \shortrightarrow 0} |u^\varepsilon(t)|^2=+\infty$ is still considered well defined; $\int{W^0(x,k,t)dk}$ may give rise to expressions that have no value, finite or otherwise}.  See \cite{Makr} for more details. This means that, if we are inetersted in an approximation of $|u^\varepsilon(x,t)|^2$, and not the WM $W^0$ itself, the WM technique does, in fact, have problems on caustics. These problems can be overcome using the SWT.

\vskip 0.25cm
This raises a more general question: can we build $\varepsilon$-dependent {\em corrections} of some sort into a phase-space asymptotic technique -- in a computable way? Recall that $W^0$ does not depend on $\varepsilon$, i.e. any WM-based technique cannot offer {\em any kind of $\varepsilon$-dependent information}. 

This is the motivation for constructing, in theorem \ref{thrm2}, the approximation $\rho_1^\varepsilon \approx \widetilde{W}^\varepsilon$, in addition to $\rho \approx W^\varepsilon$. From a mathematical point of view, $\rho(t)=W^\varepsilon(\phi_t(x,k),t=0)$ is a good approximate solution. From a numerical point of view however, the oscillations that $W^\varepsilon(x,k,t=0)$ in general contains could make even its propagation along characteristics as difficult as the direct solution of the original problem (\ref{eq1}) -- keep in mind that this is in twice the space dimensions than (\ref{eq1}). On the other hand, the oscillations of $\widetilde{W}^\varepsilon_0$ are controlled by the smoothing, making the propagation $\rho_1(t)=\widetilde{W}^\varepsilon(\widetilde{\phi}_t(x,k),t=0)$ significantly easier. (The fact the the flow comes from a smoothed potential also can simplify a number of things, e.g. in setting up a solver for the trajectories). So finally smoothing gives a {\em practical}  way to incorporate $\varepsilon$-dependent information in a phase-space asymptotic technique. 

 One should keep in mind that in general the complexity of simulating semiclassical problems explodes as $\varepsilon \shortrightarrow 0$, so apart from documenting the convergence and limitations of any mumerical scheme, a key question is to understand the rate at which complexity (computation time, memory used) grow in $\varepsilon$. See \cite{a1,a2} for more details on the use of the SWT in semiclassical simuation, including in particular a study of computational behavior in $1$-dimensional problems, of the approximation of observables (with main example $|u^\varepsilon(x,t)|^2$), and examples of calibrating the smoothing  with respect to intrinsic scales of the initial data.

\vskip 0.2cm
Finally, the fact that all the errors in this work are computed explicitly (in contrast e.g. with weak-$*$ results, where the errors depend, in a non-explicit way, on the test function) can be very helpful in setting up numerical simulations.

\subsubsection{Further work}
\label{ss003}

To conclude this long introduction, some words are in order for the more indirect motivations and implications. It has become more and more clear that there are problems where WMs, in their more straightfoward application, fundamentally fail. (For example, when there can be constructed two ($\varepsilon$-families of) problems whose WMs are identical at one point in time, but not at a later one). 

Examples include  nonlinear problems, systems with crossing eigenvalues etc -- see \cite{cfk} for a survey and a plethora of references. Another example, triggered by low smoothness in the potential, is presented already in \cite{LP}, Rem. IV 3. Loosely speaking, in all these cases the information lost in the limit to $W^0$ is ``too much''; it is worth wondering if some modified, augmented version can still be made to work. In some cases this has already been done; see e.g. \cite{cfk}. 


 The following is a version of theorem \ref{thrm2} under weaker assumptions (one can easily check that the proof still applies):

\begin{theorem} \label{thrmnext} Assume that (A1) holds for some $\theta \in (0,+\infty) \setminus \mathbb{N}$ (i.e. $\theta$ possibly smaller than $1$), and in addition that
there is a constant $M_1$ such that for all $t \in [0,T]$ 
\begin{equation}
\label{eqwass3344}
\sum\limits_{i=1}^n || \partial_{k_i}{W}_1^\varepsilon(t) ||_{L^2} \leqslant M_1 =o\left({ (\varepsilon^{-\frac{\theta}2}+\varepsilon^{-\frac{1}2}) ||\widetilde{W}^\varepsilon(t)||_{L^2} }\right).
\end{equation}
See equation (\ref{SWTeq1}) and thereafter for the precise definition of $W^\varepsilon_1(t)$.

Then there is an $O(1)$ explicit constant $C_0$ such that  for  $t\in [0,T]$ 
\begin{equation}
\label{equjmn1W}
|| \rho_1^\varepsilon(t)-\widetilde{W}^\varepsilon(t) ||_{L^2} \leqslant  T \,\, M_1 \,\,  C_0 \left({ \, \varepsilon^{\frac{\theta}2}  + \sqrt{\varepsilon} \,}\right),
\end{equation}
and
\begin{equation}
\label{eqaqjmn1W}
\frac{
|| \rho_1^\varepsilon(t)-\widetilde{W}^\varepsilon(t) ||_{L^2} }{||\widetilde{W}^\varepsilon(t) ||_{L^2} } = o(1).
\end{equation}
\end{theorem}

By modifying this kind of result (with nontrivial differences), we are able to prove a result applicable to problems with
potentials satisfying (A1), (A2) for some appropriate $\theta \in (0,1)$. It must be noted that $\theta < 1$ vs. $\theta > 1$ is a nontrivial dichotomy: $\theta \in (0,1)$ includes $C^{1} \setminus C^{1,1}$ potentials, for which there are some results, see \cite{LP}, but they show convergence to an {\em ill-posed} semiclassical limit. In other words, there are situations with $\theta <1$ where we can show that (an appropriate version of) condition (\ref{eqwass3344}) holds, which finally allows the regularization of the known, but ill-posed, limit problem [work in preparation]. A lot of the ideas and tools that had to be developed for that are in fact in this paper. 
 
 \vskip 0.2cm
 That is, apart from an end in themselves, the results shown here (including in particular the regularity results, section \ref{qrr}) are necessary in studying other, more challenging problems as well.

\vskip 1cm
The rest of the paper is organised as follows; below we define some notation we will use; in section \ref{qrr} we will present a family of results concerning the regularity of quantum (Wigner) and classical (Liouville) equations, 
trying to be as thorough as possible. The proofs can be found in section \ref{proofqrr}. 
Section \ref{proofscr} is devoted to the proofs of the theorems \ref{thrm2}, \ref{thrm23} and \ref{thrm24}. The corollaries stated earlier are straightforward applcations of these theorems, and theorem \ref{thrmnext} follows readily by retracing the proof of theorem \ref{thrm2}.

\subsection{Notations}

The Fourier transform is defined as
\begin{equation}
\widehat{f}(k)=\mathcal{F}_{x \shortrightarrow k} \left[ {f(x)} \right]
=\int\limits_{x\in \mathbb{R}^n} {e^{-2\pi i k x} f(x)dx}. 
\end{equation}
Inversion is given by
\begin{eqnarray}
\check{f}(k)=\mathcal{F}^{-1}_{x \shortrightarrow k}\left[f(x)\right]
=\int\limits_{x\in \mathbb{R}^n} {e^{2\pi ikx}f(x)dx}, \\
\mathcal{F}_{b \shortrightarrow x}\left[ \mathcal{F}^{-1}_{a \shortrightarrow b}\left[ f(a) \right] \right]=
\mathcal{F}^{-1}_{b \shortrightarrow x}\left[ \mathcal{F}_{a \shortrightarrow b}\left[ f(a) \right] \right]=
f(x).
\end{eqnarray}

\vskip 0.25cm
The Sobolev norm of order $m$ on phase-space will be defined as follows:
\begin{equation}
||f||_{H^m(\mathbb{R}^{2n})}= \sum\limits_{|a|+|b| \leqslant m} ||\partial_x^a \partial_k^b f ||_{L^2(\mathbb{R}^{2n})},
\end{equation}
where of course $a$ and $b$ are multi-indices of length $n$ each. More generally we have
\begin{equation}
||f||_{H^m(\mathbb{R}^{d})}= \sum\limits_{|a| \leqslant m} ||\partial_x^a f ||_{L^2(\mathbb{R}^{d})}.
\end{equation}
$H^{-m}$ will be the dual of $H^m$. This definition is equivalent to the more usual ones, e.g. based on $||\widehat{f}(x)(1+|x|)^m||_{L^2}$, or $ \sqrt{\sum\limits_{|a| \leqslant m} ||\partial_x^a  f ||^2_{L^2(\mathbb{R}^{2n})}}$.

\vskip 0.25cm
A more general family of Sobolev-type spaces that we will use are $\mathcal{X}^{m,p}$, with norm defined as
\begin{equation}
\label{eqaloijk}
|| f(x,k)||_{\mathcal{X}^{m,p}}= || (|X|^2+|K|^2)^{\frac{m}{2}} \widehat{f}(X,K) ||_{L^p}.
\end{equation}

\vskip 0.25cm
Moreover,
denote for future reference
\[
||\widehat{f}||_{\mathcal{F}H^s}:=||f||_{H^s}.
\]

\vskip 0.25cm
Finally, we must mention an abuse of notation that we will do: sometimes we will suppress the $\varepsilon$ dependence of certain functions for economy, using e.g. $W(t)$ instead of $W^\varepsilon(t)$, $\rho_1(t)$ instead of $\rho_1^\varepsilon(t)$ etc.

\section{Quantum and classical regularity results}\label{qrr}

The following theorem is a direct consequence of Theorem 2.1 of \cite{Mark}:
\begin{theorem}[$L^2$ regularity of the Wigner equation] \label{thrmWignL2} If the Schr\"odinger operator $-\frac{1}2\Delta +V(x)$ is essentially self-adjoint on $L^2(\mathbb{R}^n)$, then the corresponding Wigner equation preserves the $L^2$ norm, i.e. for $W^\varepsilon_0 \in L^2(\mathbb{R}^{2n})$ there is a unique solution of
\begin{equation}
\label{WignerEq}
\begin{array}{c}
\partial_t W^\varepsilon + 2\pi k \cdot \partial_x W^\varepsilon +\frac{2}{\varepsilon} Re \left[{ i\int{e^{2\pi i Sx}\hat{V}(S)W^\varepsilon(x,k-\frac{\varepsilon S}{2})dS} }\right]=0, \\
W^\varepsilon(t=0)=W^\varepsilon_0,
\end{array}
\end{equation}
and $||W^\varepsilon(t)||_{L^2(\mathbb{R}^{2n})}=||W^\varepsilon_0||_{L^2(\mathbb{R}^{2n})}$ $\forall t \in \mathbb{R}$.
\end{theorem}

\vskip 0.5cm
The following estimation is quoted from \cite{Pul}, where it appears as Theorem 4.1 (with a different but equivalent definition for $H^m$):

\noindent {\bf Statement ($H^m$ regularity of the Wigner equation)} {\em
\label{thrm1}
For each $m \in \mathbb{N}$ denote
\[
C_m=\int\limits_{\mathbb{R}^n} {|\hat{V}(S)|\,\,|S|^m dX}.
\]
Denote by $W^\varepsilon(t)$ the solution of the Wigner equation (\ref{WignerEq}), with initial data $W^\varepsilon(t=0)=W^\varepsilon_0$. Then, if $C_m<\infty$, there is a constant $C$ independent of $\varepsilon$ such that
\begin{equation}
||W^\varepsilon(t)||_{H^m} \leqslant e^{C t} ||W^\varepsilon_0||_{H^m}.
\end{equation} }

\vskip 0.8cm
{ However it contains a typo. We prove a weaker result:}

\begin{theorem}[$H^m$ regularity of the Wigner equation]\label{thrm1corr}
Denote by $W^\varepsilon(t)$ the solution of the Wigner equation (\ref{WignerEq}), with initial data $W^\varepsilon(t=0)=W^\varepsilon_0$. Then, for each $m \in \mathbb{N}$, the following estimate holds for time-scales $t\in [0,T]$, $T=O(1)$: if
\[
\int\limits_{\mathbb{R}^n} {|\hat{V}(S)|\,\, (|S|+|S|^{m+1}) dS}<\infty,
\]
then there is a constant $D=D(m,n)$ such that
\begin{equation}
\label{eqyui1}
||W^\varepsilon(t)||_{H^m} \leqslant e^{ t D(m,n)} 
||W^\varepsilon_0||_{H^m}.
\end{equation}
\end{theorem}

\noindent {\bf Remark: } See equation (\ref{eqapp001mnj}) for an estimate of the constant $D(m,n)$.

\vskip 0.5cm
The same strategy as in the proof of theorem \ref{thrm1corr} works for the Liouville equation as well:
\begin{theorem}[$H^m$ regularity for the Liouville equation]\label{liouv2}
Denote by $\rho(t)$ the solution of the Liouville equation
\begin{equation}
\label{eqllmmtrd}
 \rho_t+2\pi k \cdot \partial_x \rho-\frac{1}{2\pi} \partial_xV(x) \cdot \partial_k \rho=0
\end{equation}
with initial data $\rho(t=0)=\rho_0$. 

Assume that
\[
||\hat{V}(S)|S|^{m+1}||_{L^1} < \infty,
\]
and that, $\forall t \in[0,T]$
\begin{equation}
\label{eq9adsahg}
\begin{array}{c}
 \left|{ \langle \partial_x^A\partial_k^B \rho(t), k\cdot \partial_x \, \partial_x^A\partial_k^B \rho(t) \rangle }\right| <+\infty, \\ { } \\

 \left|{ \langle \partial_x^A\partial_k^B \rho(t), \partial_x V \cdot \partial_k \, \partial_x^A\partial_k^B \rho(t) \rangle }\right| <+\infty.
\end{array}
\end{equation}

Then the following estimate holds for $t\in [0,T]$: there is a constant $C>0$  such that
 \begin{equation}
||\rho(t)||_{H^m} \leqslant e^{ t C  }  ||\rho_0||_{H^m}
\end{equation}

\vskip 0.4cm
A sufficient condition for equation (\ref{eq9adsahg}) to hold is the following:

\noindent The  flow associated with (\ref{eqllmmtrd}) will be assumed to be complete, i.e. any bounded ball is mapped within some bounded ball for all times (trajectories do not reach infinity in finite time).
Moreover that $\rho_0$ is of compact support and,  for every multi-index $a$ with $|a|\leqslant m+2$
\begin{equation}
\label{eq8uuhvj}
\partial_x^a V(x) \in L^\infty_{loc}(\mathbb{R}^n).
\end{equation}

Another sufficient condition is that for every multi-index $a$ with $|a|\leqslant m+2$
\begin{equation}
\label{eq8uuhvj1}
\partial_x^a V(x) \in L^\infty(\mathbb{R}^n).
\end{equation}
\end{theorem}

\noindent {\bf Remark: } 
Our sufficient condition involves one more derivative in the potential (up to order $m+2$), but the estimate itself involves only derivatives of order $m+1$. This will be important for example in situations where the derivatives exist but become asymptotically large (e.g. if $V$ was a mollified version of a non-smooth function).

\vskip 0.5cm
There are some more regularity results for the Wigner equation in section \ref{proofthrmcorr2}. We don't discuss theme explicitly because we don't use them directly in the proof of our main results here. However it must be noted that in other situations they might be very useful.

 \section{Proof of regularity results}\label{proofqrr}

\subsection{Proof of Theorem \ref{thrm1corr}}

The Fourier transform of the Wigner equation (\ref{WignerEq}) is
\begin{equation}
\label{eqw1}
\hat{W}_t - 2\pi X \cdot \partial_K \hat{W} +2\int{ \hat{V}(S) \hat{W}(X-S,K) \frac{sin(\pi \varepsilon S\cdot K)}{\varepsilon} dS}=0
\end{equation}
For future reference denote by $\hat{U}(t)$ the propagator of this equation. Theorem \ref{thrmWignL2} essentially tells us that, under the self-adjointness assumption,
\begin{equation}
\forall t \in \mathbb{R} \,\,\,\,\,\,\,\,\,\,\,\,\,\,  ||\hat{U}(t)||_{L^2(\mathbb{R}^{2n}) \shortrightarrow L^2(\mathbb{R}^{2n})}=1.
\end{equation} 

Denote
\begin{equation}
v_{A,B}(X,K)=X^A K^B \hat{W}(X,K)
\end{equation}
in the usual multi-index notation.
By elementary computations, it follows that
\begin{equation}
\label{eq2}
\begin{array}{c}
X^A K^B X\cdot \partial_{K} \hat{W} = X\cdot \partial_{K} \left({v_{A,B} }\right) - \sum\limits_{B_j>0} B_j v_{A+e_j,B-e_j}
\end{array}
\end{equation}
and
\begin{equation}
\label{eq3}
\begin{array}{c}
X^AK^B \int{ \hat{V}(S) \hat{W}(X-S,K) \frac{sin(\pi \varepsilon S\cdot K)}{\varepsilon} dS}= \\ { } \\
= \int{ \hat{V}(S) v_{A,B}(X-S,K) \frac{sin(\pi \varepsilon S\cdot K)}{\varepsilon} dS}+ \\ { } \\
+  \sum\limits_{0<l\leqslant A} { \prod\limits_{i=1}^n \binom{A_i}{l_i}
  \int{ \hat{V}(S)  S^{l}    v_{A-l,B}(X-S,K) \frac{sin(\pi \varepsilon S\cdot K)}{\varepsilon} dS} }.
\end{array}
\end{equation}

To see that, it suffices to observe that
\[
\begin{array}{c}
X^A=(S+(X-S))^A=\sum\limits_{l=0}^A \binom{A}{l} S^l (X-S)^{A-l} \,\, \Rightarrow 
\\ { } \\
 \Rightarrow\,\, X^A=(X-S)^A + \sum\limits_{l\in \mathbb{L}(A)} \binom{A}{l} S^l (X-S)^{A-l} 
\end{array}
\]
where for brevity we use  the notation
\[
\binom{A}{l}
:=\prod\limits_{i=1}^n \binom{A_i}{l_i},
\]
and
\begin{equation}
\label{eqppaammzz}
\mathbb{L}(A) := \{ { l \in \mathbb{N}^n |  l \leqslant A, \,\, |l|>0 \} }\}.
\end{equation}
(Remark on notation: $l\leqslant A \Leftrightarrow \forall i \in \{1,...,n\} : l_i \leqslant A_i$, and $|l|=\sum\limits_{i=1}^n l_i$). We will use these notations freely in the sequel.

Two things should be noted: obviously if $A_i=0$ then there is no contribution from the $i$ coordinate; and $l \in \mathbb{L}(A) \Rightarrow 0<|l|\leqslant |A|$.

So multiplying equation (\ref{eqw1}) with $X^A K^B$ yields
\begin{equation}
\label{eq3aa}
\begin{array}{c}
\frac{d}{dt}v_{A,B} - 2\pi X \cdot \partial_K v_{A,B} +2\int{ \hat{V}(S) v_{A,B}(X-S,K) \frac{sin(\pi \varepsilon S\cdot K)}{\varepsilon} dS}=\\ { } \\

=-2\pi \sum\limits_{B_j>0} B_j v_{A+e_j,B-e_j}-\\ { } \\

-2\sum\limits_{l \in \mathbb{L}(A) } \binom{A}{l} \int{ \hat{V}(S)S^{l} v_{A-l,B}(X-S,K) \frac{sin(\pi \varepsilon S\cdot K)}{\varepsilon} dS},
\end{array}
\end{equation}
which, for any $m \in \mathbb{N}$, is a closed system for $v_{A,B}$ with $|A+B|\leqslant m$. (For $A=0$, we just forget the last term).

Moreover, observe that
\begin{equation}
\label{eq8ui9oklo0}
\frac{d}{dt}||v_{A,B}(t)||_{L^2}^2=2||v_{A,B}(t)||_{L^2} \,\, \frac{d}{dt}|| v_{A,B}(t)||_{L^2}
\end{equation}
while at the same time, using equation (\ref{eq3aa}),
\begin{equation}
\label{eq8ui111}
\begin{array}{c}
\frac{d}{dt}||v_{A,B}(t)||_{L^2}^2= \langle{v_{A,B},\partial_t v_{A,B} }\rangle+ \langle{\partial_t v_{A,B},v_{A,B} }\rangle=\\ { } \\

= 
2 Re \left[{
\langle 
v_{A,B}(t),  2\pi X \cdot \partial_K v_{A,B} \rangle - }\right. \\ { } \\
-\langle v_{A,B}(t), 2\int{ \hat{V}(S) v_{A,B}(X-S,K) \frac{sin(\pi \varepsilon S\cdot K)}{\varepsilon} dS}\rangle+ \\ { } \\
+\langle v_{A,B}, 2\pi \sum\limits_{B_j>0} B_j v_{A+e_j,B-e_j} \rangle -  \\ { } \\

\left.{
-\langle v_{A,B}, \sum\limits_{l \in \mathbb{L}(A) } \binom{A}{l} 2\int{ \hat{V}(S)S^{l} v_{A-l,B}(X-S,K) \frac{sin(\pi \varepsilon S\cdot K)}{\varepsilon} dS} \rangle  }\right]
\end{array}
\end{equation}

Assuming for a moment that the integrals exist, we observe that
\begin{equation}
\label{eq9iok341}
Re \left[{ \langle v_{A,B}(t),  2\pi X \cdot \partial_K v_{A,B} \rangle }\right]=0,
\end{equation}
and
\begin{equation}
\label{eq9iok34}
Re \left[{ \langle v_{A,B}(t),  2\int{ \hat{V}(S) v_{A,B}(X-S,K) \frac{sin(\pi \varepsilon S\cdot K)}{\varepsilon} dS} \rangle }\right]=0
\end{equation}
by anti-symmetry.

Equation (\ref{eq9iok341}) is completely obvious; for (\ref{eq9iok34}) observe that
\[
\begin{array}{c}
Re \left[{ \int{ \overline{f}(X,K)  2\int{ \hat{V}(S) f(X-S,K) \frac{sin(\pi \varepsilon S\cdot K)}{\varepsilon} dS} dXdK} }\right]=\\ { } \\

=\int{ \overline{f}(X,K)  { \hat{V}(S) f(X-S,K) \frac{sin(\pi \varepsilon S\cdot K)}{\varepsilon} dS}dXdK }+\\
+\int{ {f}(X,K)  \overline{{ \hat{V}(S) f(X-S,K) \frac{sin(\pi \varepsilon S\cdot K)}{\varepsilon} dS}dXdK} }=\\ { } \\

=\int{ \overline{f}(X,K)  { \hat{V}(S) f(X-S,K) \frac{sin(\pi \varepsilon S\cdot K)}{\varepsilon} dS}dXdK }+\\
+\int{ {f}(X,K)  { \hat{V}(-S) \overline{f}(X-S,K) \frac{sin(\pi \varepsilon S\cdot K)}{\varepsilon} dS}dXdK} =\\ { } \\

=\int{ \overline{f}(X,K)  { \hat{V}(S) f(X-S,K) \frac{sin(\pi \varepsilon S\cdot K)}{\varepsilon} dS}dXdK }-\\
-\int{ {f}(X,K)  { \hat{V}(S) \overline{f}(X+S,K) \frac{sin(\pi \varepsilon S\cdot K)}{\varepsilon} dS}dXdK} =\\ { } \\

=\int{ \overline{f}(X,K)  { \hat{V}(S) f(X-S,K) \frac{sin(\pi \varepsilon S\cdot K)}{\varepsilon} dS}dXdK }-\\
-\int{ {f}(X-S,K)  { \hat{V}(S) \overline{f}(X,K) \frac{sin(\pi \varepsilon S\cdot K)}{\varepsilon} dS}dXdK} =0
\end{array}
\]

Of course for these considerations to be valid the integrals must exist; they do, thanks to the following lemma:

\begin{lemma}[$H^m_\varepsilon$ regularity of the Wigner equation]\label{thrm7yu8i9o}
Denote by $W(t)$ the solution of the Wigner equation like earlier.
Assume that
\begin{equation}
\label{eq9876jhgf}
\int\limits_{\mathbb{R}^n} {|\hat{V}(S)|\,\, (|S|+|S|^{m}) dS}<\infty, \,\,\,\, ||W^\varepsilon_0||_{H^m}< \varepsilon^{-c}
\end{equation}
for some $c>0$.

Then, for each $m \in \mathbb{N}$, $T>0$, there exist constants $C,\mu>0$ such that
 \begin{equation}
 \label{equ8i1lm}
||W^\varepsilon(t)||_{H^m} \leqslant C \varepsilon^{-\mu}, \\ { } \\
\end{equation}
and
 \begin{equation}
\forall |A+B| \leqslant m-1, \,\, i\in \{ 1,...,n\}, \,\,\,\, ||k_i \partial_x^A\partial_k^B W(t)||_{L^2} \leqslant C \varepsilon^{-\mu}.
\end{equation}
\end{lemma}

(For our purposes it is obvious the explicit computation of $C$, $\mu$ is not important -- it is possible, as one can readily see in the proof of lemma \ref{thrm7yu8i9o}, in section \ref{proofthrmcorr2}. Observe moreover that the requirement $\widehat{V}(S)|S| \in L^1$, which allows for weaker possible blowup at $0$ than $\widehat{V}(S)|S|^{m+1} \in L^1$, originally comes from lemma \ref{thrm7yu8i9o}).

\vskip 0.7cm
Now, by combining equations (\ref{eq8ui9oklo0}) and (\ref{eq8ui111}) and discarding the terms (\ref{eq9iok341}), (\ref{eq9iok34}) that vanish identically, it follows that\footnote{We also used the obvious estimate $\left|{\hat{V}(S)S^{l} v_{A-l,B}(X-S,K) \frac{sin(\pi \varepsilon S\cdot K)}{\varepsilon}}\right| \leqslant \left|{\hat{V}(S)S^{l} v_{A-l,B}(X-S,K) \pi S \cdot K}\right|$ in equation (\ref{eq8ui111}).} 

\begin{equation}
\label{wgradInteq2}
\begin{array}{c}
\frac{d}{dt}||v_{A,B}(t)||_{L^2} \leqslant    \pi \sum\limits_{B_j>0} B_j ||v_{A+e_j,B-e_j}(t)||_{L^2} +\\ { } \\

+\pi  \sum\limits_{l \in \mathbb{L}(A) } \binom{A}{l} { \sum\limits_{j=1}^n ||\hat{V}(S)S^{l+e_j}||_{L^1} ||v_{A-l,B+e_j}(t)||_{L^2} } .

\end{array}
\end{equation}

Now we sum up all the equations for $|A|+|B|\leqslant m$ to get

\begin{equation}
 \begin{array}{c}
\label{eqddt}
  \frac{d}{dt} \sum\limits_{|A+B|\leqslant m} {||v_{A,B}(t)||_{L^2}} \leqslant  \\ { } \\

\leqslant n^{m+1} \pi m!  \cdot max \{ 1, ||\hat{V}(S)|S|^{m+1}||_{L^1} \}  \, \sum\limits_{|A+B|\leqslant m} {||v_{A,B}(t)||_{L^2}}.
 \end{array}
\end{equation}

With the help of Gronwall's lemma the result follows. For brevity, in the sequel we will denote
\begin{equation}
\label{eqapp001mnj}
D(m,n):=(1+n^{m+1}  m!)\pi  \cdot max \{ 1, ||\hat{V}(S)(|S|+|S|^{m+1})||_{L^1} \}.
\end{equation}

\noindent {\bf Remark:}
For $m=1$ it is very easy to see that we can have a slightly better constant, namely:
\begin{equation}
 \begin{array}{c}
\label{eqddt00}
  \frac{d}{dt} \sum\limits_{|A+B|\leqslant 1} {||v_{A,B}(t)||_{L^2}} \leqslant  \\ { } \\

\leqslant n \pi \, max \{ 1, ||\hat{V}(S)|S|^{2}||_{L^1} \}  \, \sum\limits_{|A+B|\leqslant 1} {||v_{A,B}(t)||_{L^2}},
 \end{array}
\end{equation}
i.e. $D(1,n)=n \pi \, max \{ 1, ||\hat{V}(S)|S|^{2}||_{L^1} \}$.

 \subsection{Proof of Theorem \ref{liouv2}}\label{proofliouv2}
 \noindent {\bf Proof:} The proof is essentially the same as for theorem \ref{thrm1corr}: we have 
\[
\pi \int{ \hat{V}(S) \rho(X-S,K) S\cdot K dS}
\]
instead of 
\[
\int{ \hat{V}(S) \hat{W}(X-S,K) \frac{sin(\pi \varepsilon S\cdot K)}{\varepsilon} dS}
\]
in equation (\ref{eqw1}), and then the proof follows along similar lines.
The
counterparts of the terms equations (\ref{eq9iok341}), (\ref{eq9iok34}) are
\begin{equation}
\label{eq9akjyhg}
\begin{array}{c}
Re \left[{ \langle \partial_x^A\partial_k^B \rho, k\cdot \partial_x \, \partial_x^A\partial_k^B \rho \rangle }\right], \\ { } \\

Re \left[{ \langle \partial_x^A\partial_k^B \rho, \partial_x V \cdot \partial_k \, \partial_x^A\partial_k^B \rho \rangle }\right].
\end{array}
\end{equation}

It is clear that, if the integrals exist, these terms vanish identically by anti-symmetry.

Now we will work out our sufficient condition for the existence of the integrals of equation (\ref{eq9akjyhg}):
it follows easily,  using the method of characteristics, that if 
\[
\mathop{sup}\limits_{
\begin{array}{c}
|a|\leqslant m+2 \\
x \in \mathbb{R}^n
\end{array}
} |\partial_x^a V(x)|<\infty,
\]
 and the flow is {\em complete} (see the statement of the theorem), then, for initial data of compact support such that
\begin{equation}
\partial_x \partial_x^A\partial_k^B \rho(0), \,\, \partial_x \partial_x^A\partial_k^B \rho(0), \,\, k \partial_x \partial_x^A\partial_k^B \rho(0) \,\, \in L^2,
\end{equation}
 we have $\forall t>0$
 \begin{equation}
\partial_x \partial_x^A\partial_k^B \rho(t), \,\, \partial_x \partial_x^A\partial_k^B \rho(t), \,\, k \partial_x \partial_x^A\partial_k^B \rho(t) \,\, \in L^2,
\end{equation}
and hence the integrals in question exist.
  
  The second sufficient condition follows clearly.
  
  The proof is complete.

 \subsection{Auxiliary lemmata}\label{proofthrmcorr2}

It is more convenient to prove lemma \ref{thrm7yu8i9o} first for $m=1$ and then for $m \in \mathbb{N}$:

\begin{lemma}[$H^1_\varepsilon$ regularity of the Wigner equation]\label{thrm1corr2}
Denote by $W^\varepsilon(t)$ the solution of the Wigner equation (\ref{WignerEq}), with initial data $W^\varepsilon(t=0)=W^\varepsilon_0$. Then the following estimate holds for $t\in [0,T]$:
if
\[
\int\limits_{\mathbb{R}^n} {|\hat{V}(S)|\,\, |S|dS}<\infty,
\]
then we have
\begin{equation}
\begin{array}{c}
||\varepsilon \partial_{x_i} W^\varepsilon(t)||_{L^2}\leqslant ||\varepsilon \partial_{x_i} W^\varepsilon_0||_{L^2}+2t || \widehat{V}(S)|S|\,||_{L^1} ||W^\varepsilon_0||_{L^2}
 \\ { } \\

||\varepsilon \partial_{k_i} W^\varepsilon(t)||_{L^2}\leqslant 
||\varepsilon \partial_{k_i} W^\varepsilon_0||_{L^2}+ \,\,\,\,\,\,\,\,\,\,\,\,\,\,\,\,\,\,\,\,\,\,\,\,\,\,\,\,\,\,\,\,\, \\
\,\,\,\,\,\,\,\,\,\,\,\,\,\,\,\,\,\,\,\,\,\,\,\,\,\,\,\,\,\,\,\,\,
+2\pi t \left({ ||\partial_{x_i} W^\varepsilon_0||_{L^2}+2\frac{t}{\varepsilon} || \widehat{V}(S)|S|\,||_{L^1} ||W^\varepsilon_0||_{L^2} }\right)
, \\ { } \\

||\varepsilon k_i  W^\varepsilon(t)||_{L^2}\leqslant ||\varepsilon k_i  W^\varepsilon_0||_{L^2}+2\pi t || \widehat{V}(S)|S|\,||_{L^1} ||W^\varepsilon_0||_{L^2}.
\end{array}
\end{equation}
\end{lemma}

\vskip 0.8cm
\noindent {\bf Proof of lemma \ref{thrm1corr2}:}
 
Denote $u_i(X,K)=X_i \widehat{W}(X,K)$, $v_i(X,K)=K_i \widehat{W}(X,K)$, $z_i=\partial_{K_i}\widehat{W}(X,K)$. Then one checks that equation (\ref{eqw1}) implies
\begin{equation}
\label{eq8uii2w}
\begin{array}{c}
\partial_t u_i - 2\pi X \cdot \partial_K u_i +2\int{ \hat{V}(S) u_i(X-S,K) \frac{sin(\pi \varepsilon S\cdot K)}{\varepsilon} dS}=\\
=-2\int{\widehat{V}(S)S_i\widehat{W}(X-S,K)\frac{sin(\pi \varepsilon S K)}{\varepsilon}dS}, \\ { } \\

\partial_t v_i - 2\pi X \cdot \partial_K v_i +2\int{ \hat{V}(S) v_i(X-S,K) \frac{sin(\pi \varepsilon S\cdot K)}{\varepsilon} dS}=\\
=2\pi u_i, \\ { } \\

\partial_t z_i - 2\pi X \cdot \partial_K z_i +2\int{ \hat{V}(S) z_i(X-S,K) \frac{sin(\pi \varepsilon S\cdot K)}{\varepsilon} dS}=\\
=-2\pi \int{ \hat{V}(S)S_i \widehat{W}(X-S,K) \frac{cos(\pi \varepsilon S\cdot K)}{\varepsilon} dS}.
\end{array}
\end{equation}

Making use of the second part of lemma \ref{lmANS} (equation (\ref{eqansb})), it follows that
an equivalent reformulation of equation (\ref{eq8uii2w}) in terms of the propagator $\widehat{U}(t)$ is
\begin{equation}
\begin{array}{c}
\label{eqokmnjhu}
u_i(t)=\widehat{U}(t)u_i(0)+2\int\limits_{\tau=0}^t {\widehat{U}(t-\tau)\int{\widehat{V}(S)S_i\widehat{W}(X-S,K,\tau)\frac{sin(\pi \varepsilon S K)}{\varepsilon}dS}d\tau}
, \\ { } \\

v_i(t)=\widehat{U}(t)v_i(0)-\int\limits_{\tau=0}^t {\widehat{U}(t-\tau) 2\pi u_i(\tau)d\tau}, \\ { } \\

z_i(t)=\widehat{U}(t)z_i(0)+2\pi\int\limits_{\tau=0}^t {\widehat{U}(t-\tau)  \int{ \hat{V}(S)S_i \widehat{W}(X-S,K,\tau) \frac{cos(\pi \varepsilon S\cdot K)}{\varepsilon} dS}d\tau}.
\end{array}
\end{equation}
Now, making use of theorem \ref{thrmWignL2}, the obvious bounds $|sin(\pi\varepsilon S K )|,cos(\pi\varepsilon S K )|\leqslant 1$),
 and lemma \ref{lemyoungineq},
we readily deduce from equation (\ref{eqokmnjhu}) that
\begin{equation}
\label{equation69}
\begin{array}{c}
||u_i(t)||_{L^2}\leqslant ||u_i(0)||_{L^2}+2\frac{t}{\varepsilon} || \widehat{V}(S)|S|\,||_{L^1} ||W^\varepsilon_0||_{L^2}
 \\ { } \\

||v_i(t)||_{L^2}\leqslant ||v_i(0)||_{L^2}+2\pi t \mathop{sup}\limits_{\tau\in(0,t)} ||u_i(\tau)||_{L^2} \leqslant\\
\leqslant ||v_i(0)||_{L^2}+2\pi t \left({ ||u_i(0)||_{L^2}+2\frac{t}{\varepsilon} || \widehat{V}(S)|S|\,||_{L^1} ||W^\varepsilon_0||_{L^2} }\right)
, \\ { } \\

||z_i(t)||_{L^2}\leqslant ||z_i(0)||_{L^2}+\frac{2\pi t}{\varepsilon} || \widehat{V}(S)|S|\,||_{L^1} ||W^\varepsilon_0||_{L^2}.
\end{array}
\end{equation}

The proof is complete.

\vskip 0.8cm
Now we are ready for the general case.

\noindent {\bf Proof of lemma \ref{thrm7yu8i9o}:}

 System (\ref{eq3aa}) can equivalently be recast as
\begin{equation}
\label{wgradInteq}
\begin{array}{c}
v_{A,B}(t)=\hat{U}(t)v_{A,B}(0) 

-\int\limits_{\tau=0}^t {  \hat{U}(t-\tau) \left({2\pi \sum\limits_{B_j>0} B_j v_{A+e_j,B-e_j}(\tau) + }\right. } \\ { } \\

+2\left.{ 
\sum\limits_{l \in \mathbb{L}(A) } \binom{A}{l} \int{ \hat{V}(S)S^{l} v_{A-l,B}(X-S,K,\tau) \frac{sin(\pi \varepsilon S\cdot K)}{\varepsilon} dS} 
}\right)d\tau.
\end{array}
\end{equation}

The point to repeating the idea of the previous proof, is to go through the multiindices $(A,B)$ so that we always have on the right-hand-side a quantity that we have already estimated. 

This is easy to check; for example let us prove the result for $m>1$ assuming it holds for $m-1$. First work with all the indices such that $|A|+|B|=m$, $|B|=0$. It is clear that, for those terms equation (\ref{wgradInteq}) is just
\begin{equation}
\begin{array}{c}
v_{A,0}(t)=\hat{U}(t)v_{A,0}(0) +\\ { } \\

+2\int\limits_{\tau=0}^t{ 
\sum\limits_{l \in \mathbb{L}(A) } \binom{A}{l} \int{ \hat{V}(S)S^{l} v_{A-l,0}(X-S,K,\tau) \frac{sin(\pi \varepsilon S\cdot K)}{\varepsilon} dS} 
d\tau}.
\end{array}
\end{equation}
But $|l|>0$ (recall the definition of $\mathbb{L}(A)$, equation (\ref{eqppaammzz})) means that all the terms of the form $v_{A-l,0}$ in the right-hand-side are of order at most $m-1$, and therefore bounded by assumption. By using once again lemma \ref{lemyoungineq}, and integrating in time, we get 
\begin{equation}
\begin{array}{c}
||v_{A,0}(t)||_{L^2} \leqslant ||v_{A,0}(t)||_{L^2} +\\
+2
\frac{t}{\varepsilon} 
\sum\limits_{l \in \mathbb{L}(A) } \binom{A}{l}

||\widehat{V}(S)S^{l}||_{L^1}
 \mathop{sup}\limits_{
\begin{footnotesize}\begin{array}{c}
\tau \in (0,t)
\end{array}
\end{footnotesize}} ||v_{A-l,0}(\tau)||_{L^2} \leqslant \\ { } \\

\leqslant ||v_{A,0}(t)||_{L^2} +\\
+2
\frac{t}{\varepsilon} 
||\widehat{V}(S)|S|^{m}||_{L^1}

\sum\limits_{l \in \mathbb{L}(A) } \binom{A}{l}

 \mathop{sup}\limits_{
\begin{footnotesize}\begin{array}{c}
\tau \in (0,t)
\end{array}
\end{footnotesize}} ||v_{A-l,0}(\tau)||_{L^2}.
\end{array}
\end{equation}
Now we can allow $|A|+|B|=m$, $|B|=1$, and observe that the rhs contains only terms of the form $|A|+|B|=m$, $|B|=0$, which we just estimated. To conclude we just proceed inductively until $|B|=m$.

The second inequality in (\ref{equ8i1lm}) follows similarly, by a straightforward adaptation of the previous proof. 

The result now follows.

\vskip 0.3cm
\begin{lemma}[$H^{-m}$ regularity of the Wigner equation] \label{lmapqolz} Under the assumptions of theorem \ref{thrm1corr},
\[
||W^\varepsilon(t)||_{H^{-m}} \leqslant e^{tD(m,n)}||W^\varepsilon_0||_{H^{-m}}.
\]
\end{lemma}

\noindent {\bf Proof of lemma \ref{lmapqolz}:} The result follows directly by noting that
\[
\mathop{sup}\limits_{
 ||g||_{H^m} \leqslant 1
} |\langle{ U(t)f,g }\rangle| =
\mathop{sup}\limits_{
 ||g||_{H^m} \leqslant 1
} |\langle{ f,U(-t)g }\rangle| \leqslant ||U(-t)||_{H^m \shortrightarrow H^m} ||f||_{H^{-m}}.
\]

\vskip 0.3cm

Finally, the following regularity result can be shown with minimal assumptions on the potential -- although apparently it cannot be generalized to higher order derivatives:

\begin{lemma} Assume there is a $C_*>0$ such that $V(x)>-C_*$, and the Schr\"odinger operator $\frac{\varepsilon^2}{2}\Delta-V(x)$ is self adjoint.

Moreover, consider a mixed state as follows\footnote{A possible example would be of the form $\sum\limits_{m\in\mathbb{N}} {\lambda_m u^\varepsilon_m(x)\overline{u}^\varepsilon_m(y)}$, with appropriate decay of $\lambda_m$, regularity of $u^\varepsilon_m$ of course.}: denote $\lambda \in \Lambda$ an appropriate index set, and
\begin{equation}
\label{equationanscs}
\rho^\varepsilon_0(x,y)=\int{ u^\varepsilon_\lambda(x)\overline{u}^\varepsilon_\lambda(y) d\mu(\lambda)},
\end{equation}
where $||u_\lambda^\varepsilon(t=0)||_{L^2}=1$, and
\[
|\langle V(x)u^\varepsilon_\lambda(t=0), u^\varepsilon_\lambda(t=0) \rangle|+
\sum\limits_{i=1}^n {\left({ ||\partial_{x_i}u_\lambda^\varepsilon(t=0)||_{L^2} + ||x_iu_\lambda^\varepsilon(t=0)||_{L^2} }\right)} \leqslant C_\lambda .
\]
Moreover, assume $(1+C_\lambda) \in L^1(d\mu(\lambda),\Lambda)$. (In particular no assumption is made for the behaviour of $C_\lambda$ in $\varepsilon$).

Naturally the relation between $\rho^\varepsilon$ and the Wigner function is as in equation (\ref{eqdefrelrowig}).

Then there are $C',c>0$ such that for all $i=1,...,n$, $t\in[0,T]$,
\[
\begin{array}{c}
||\partial_{x_i}W^\varepsilon(t)||_{L^2} \leqslant  C' \varepsilon^{-c}, \\ { } \\

||\partial_{k_i}W^\varepsilon(t)||_{L^2} \leqslant  C'(1+T) \varepsilon^{-c}.
\end{array}
\]
\end{lemma}

\noindent {\bf Proof:} First of all, if $u^\varepsilon_\lambda(x,t)$ is the solution of (\ref{eq1}) with initial data $u^\varepsilon_\lambda(x)$, it is straightforward to check that
\[
\rho^\varepsilon(x,y,t)=\int{ u^\varepsilon_\lambda(x,t)\overline{u}^\varepsilon_\lambda(y,t) d\mu(\lambda)},
\]
i.e. the measure $d\mu(\lambda)$ doesn't depend on time.

Now recall the energy conservation for equation (\ref{eq1}): 
\begin{equation}
\begin{array}{c}
\langle \frac{\varepsilon^2}{2}\Delta u^\varepsilon(t)-V(x)u^\varepsilon(t),u^\varepsilon(t)\rangle = \langle \frac{\varepsilon^2}{2}\Delta u_0^\varepsilon-V(x)u_0^\varepsilon,u_0^\varepsilon\rangle \,\, \Rightarrow \\ {  } \\

\Rightarrow ||\nabla u^\varepsilon(t)||_{L^2}^2=\langle -\Delta u^\varepsilon(t),u^\varepsilon(t) \rangle=2\frac{ \langle- \frac{\varepsilon^2}{2}\Delta u_0^\varepsilon+V(x)u_0^\varepsilon,u_0^\varepsilon\rangle -\langle Vu^\varepsilon(t),u^\varepsilon(t) \rangle }{\varepsilon^2} \Rightarrow \\ { } \\

\Rightarrow  ||\nabla u^\varepsilon(t)||_{L^2} \leqslant 2\frac{\sqrt{ \langle -\frac{\varepsilon^2}{2}\Delta u_0^\varepsilon+V(x)u_0^\varepsilon,u_0^\varepsilon\rangle +C_*||u^\varepsilon_0||_{L^2} }}{\varepsilon}.
\end{array}
\end{equation}

Applying that to $u^\varepsilon_\lambda$, it follows that there is a $B>0$ independent of $\lambda,\varepsilon$ so that
\begin{equation}
\label{eqencoscff}
||u^\varepsilon_\lambda(t)||_{H^1} \leqslant B\frac{C_\lambda+C_*  }{\varepsilon}
\end{equation}
for all $t \in \mathbb{R}$. Now observe that (in the notations of lemma \ref{lemalan}) in general
\[
\partial_{x_i} W^\varepsilon = Q\left[{(\partial_{x_i}+\partial_{y_i}) \rho^\varepsilon }\right]
\]
(i.e. not only for pure states). Therefore (using equation (\ref{eqwwwnnnn22}))
\[
\begin{array}{c}
|| \partial_{x_i} W^\varepsilon(t)||_{L^2} = \varepsilon^{-\frac{n}{2}} 2\pi||(\partial_{x_i}+\partial_{y_i})\rho^\varepsilon(t)||_{L^2}\leqslant \\ { } \\

\leqslant\varepsilon^{-\frac{n}{2}} 2\pi \int{||(\partial_{x_i}+\partial_{y_i}) u^\varepsilon_\lambda(x,t)\overline{u}^\varepsilon_\lambda(y,t) ||_{L^2} d\mu(\lambda)} \leqslant  \\ { } \\

\leqslant \varepsilon^{-\frac{n}{2}-1} 4\pi \int{(C_\lambda+C_*)d\mu(\lambda)} \leqslant B' \varepsilon^{-c}
\end{array}
\]
for some $B',c>0$. (We used the joint constraint on $C_\lambda, d\mu(\lambda)$ to interchange norms and derivatives with the integration).

Now for the $\partial_{k_i}$ derivatives, recall 
\[ 
||\partial_{k_i}W^\varepsilon(t)||_{L^2}\leqslant ||\partial_{k_i}W^\varepsilon(0)||_{L^2}+2\pi t \mathop{sup}\limits_{\tau\in(0,t)} ||\partial_{x_i}W^\varepsilon(\tau)||_{L^2} ;
\]
see equation (\ref{equation69}) for the derivation. The proof is complete.

\section{Proof of the main  results}\label{proofscr}
 
\subsection{Proof of Theorem \ref{thrm2}} 
 
 The proof will be broken down to several steps. Denote be $\widetilde{W}_1^\varepsilon$ the the SWT corresponding to the problem with the potential $V_1$, i.e. the solution to
\begin{equation}
\label{SWTeq1}
\begin{array}{l}
\partial_t \widetilde{W}_1^\varepsilon + \left({ 2\pi k \cdot \partial_x +\frac{\varepsilon\sigma_x^2}{2}\partial_x\cdot\partial_k }\right) \widetilde{W}_1^\varepsilon+ \\ 

\,\,\,\,\,\,\,\,\,\,\,\,\,\,\,\,\,
+\frac{2}{\varepsilon} Re \left[{ 
i\int{e^{2\pi i Sx-\frac{\varepsilon\pi}{2}\sigma_x^2S^2}\hat{V}_1(S)\widetilde{W}_1^\varepsilon(x+\frac{i\varepsilon\sigma_x^2S}{2},k-\frac{\varepsilon S}{2})dS}  }\right]=0, \\ { } \\

\,\,\,\,\,\,\,\,\,\,\,\,\,\,\,\,\,\,\,\,\,\,\,\,\,\,\,\,\,\,\,\,\,\,\,\,\,\,\,\,\,\,\,\,\,\,\,\,\,\,\,\,\,\,\,\,\,\,\,\,\,\,\,\,\,\,
\widetilde{W}_1(x,k,0)=\widetilde{W}^\varepsilon_0.
\end{array}
\end{equation}
See equation (\ref{eqV1}) for the definition of $V_1$.
Naturally the corresponding Wigner function will be denoted by $W_1^\varepsilon(t)=\Phi^{-1} \widetilde{W}^\varepsilon_1(t)$.
Moreover denote by $\widetilde{W}^\varepsilon(t)$ the exact SWT satisfying equation (\ref{SWTeq}), $\widetilde{W}^\varepsilon(t)=\Phi W^\varepsilon(t)$. Finally, recall the definition of $\rho_1^\varepsilon$, equation (\ref{SmLiouvEqa}), and denote for further use
\begin{equation}
\label{eqdefm1}
M_1=\sum\limits_{i=1}^n \mathop{sup}\limits_{t\in [0,T]} ||\,\,|2\pi K_i| \widehat{W}_1(X,K,t)||_{L^2}.
\end{equation}

$\bullet$ \,\, {\bf Proof of equations (\ref{equjmn1}),  (\ref{equjmn2}):}

\vskip 0.25cm
There is a constant $C$, depending only on $n,\sigma_x,\sigma_k$, such that
\begin{equation}
\label{eqPart1}
||  \widetilde{W}_1^\varepsilon   -\rho_1^\varepsilon ||_{L^2} \leqslant CT(\varepsilon^{\frac{\theta}2}M_1+\sqrt{\varepsilon} M_1),
\end{equation}
\begin{equation}
\label{eqPart2}
||  \widetilde{W}_1^\varepsilon   -\widetilde{W}^\varepsilon||_{L^2} \leqslant CT(\varepsilon^{\frac{\theta}2}M_1+\varepsilon M_1),
\end{equation}
\begin{equation}
\label{eqPart3}
||  \widetilde{W}^\varepsilon(t)   -W^\varepsilon(t)||_{L^2} \leqslant C e^{tD(1,n)}\sqrt{\varepsilon} ||W_0^\varepsilon||_{H^1} , 
\end{equation}
and 
\begin{equation}
\label{eqPart3a}
|| \rho_1^\varepsilon(t)- \rho(t)||_{L^2} \leqslant  CT (\varepsilon+\varepsilon^{\frac{\theta}2}) e^{TD(1,n)}  ||W_0^\varepsilon||_{H^1}.
\end{equation}
(In the process of this last proof we will have to use that $\theta>1$). Recall that $D(1,n)$ was estimated in equation (\ref{eqddt00}). 

The  conclusion then follows by estimating
\begin{equation}
\label{eq00iibbvv}
M_1 \leqslant ||W_1^\varepsilon(t)||_{H^1} \leqslant e^{tD(1,n)} ||W^\varepsilon_0||_{H^1},
\end{equation}
where theorem \ref{thrm1corr} was used for the last inequality -- making use of A1(r=0).

$\bullet$ \,\, {\bf Proof of equations (\ref{eqaqjmn1}), (\ref{eqaqjmn2}):}

It suffices to combine equations (\ref{equjmn1}),  (\ref{equjmn2}), theorem \ref{thrmWignL2}, and to use once again equation (\ref{eqPart3}), to get 
\begin{equation}
||\widetilde{W}^\varepsilon(t)||_{L^2}=||W^\varepsilon_0||_{L^2} + \sqrt{\varepsilon} \, \zeta  e^{TD(1,n)} ||W^\varepsilon_0||_{H^1} 
\end{equation}
with $|\zeta| \leqslant C$.

\vskip 0.5cm
Now we proceed to the proofs of the building blocks:

\noindent {\bf Proof of equation (\ref{eqPart1}):} 

One can Taylor expand $\widetilde{W}_1^\varepsilon$ into a power series \footnote{because $\widetilde{W}^\varepsilon$ is an entire analytic function; see \cite{AMP} for proof and more details. This is in contrast to the {\em formal expansion} often used for the Wigner transform.} to recast equation (\ref{SWTeq1}) as
\begin{equation}
\begin{array}{l}
\partial_t \widetilde{W}_1^\varepsilon + 2\pi k \cdot \partial_x  \widetilde{W}_1^\varepsilon -
\frac{1}{2\pi} \partial_x \widetilde{V}_1(x) \cdot \partial_k \widetilde{W}_1^\varepsilon=F \widetilde{W}_1^\varepsilon
\end{array}
\end{equation}
where
\begin{equation}
\label{eqF}
\begin{array}{l}
F=-\frac{\varepsilon\sigma_x^2}{2}\partial_x\cdot\partial_k  -
\\ { } \\
-2\sum\limits_{m=2}^{\infty} 
{ 
\frac{\varepsilon^{m-1}}{(4\pi)^m} 
\sum\limits_{(m-l)mod2=1} {
 i^{l-m+1} \sigma_x^{2l} (-1)^{m-l} 
  \sum\limits_{
\begin{scriptsize}\begin{array}{c}
|A|=l \\
|B|=m-l
\end{array}
\end{scriptsize}} 
 { \frac{ \partial_x^{A+B} \widetilde{V}_1(x) }{A!B!} \partial_x^A\partial_k^B  }  
}}
\end{array}
\end{equation}
and $A,B \in (\mathbb{N} \cup \{0 \})^n$ are multi-indices. See  lemma \ref{lemJusInter} for the full justification of this step (i.e. of interchanging the summation of the Taylor expansion and the $dS$ integration).

It is clear now that the SWT equation (\ref{SWTeq1}) can be seen as a perturbation of the modified Liouville equation (\ref{SmLiouvEqa}). We will use lemma \ref{lmANS}; however before that some remarks are in order: 

Denote $\widetilde{E}_1(t)$ the propagator of equation (\ref{SmLiouvEqa}), i.e. $\rho_1^\varepsilon(t)=\widetilde{E}_1(t)\widetilde{W}^\varepsilon_0$, and $\widetilde{U}_1(t)$ the propagator of equation (\ref{SWTeq1}), $\widetilde{W}_1^\varepsilon(t)=\widetilde{U}_1(t)\widetilde{W}^\varepsilon_0$. 
The Liouville equation (\ref{SmLiouvEqa}) can be solved with the method of characteristics for any $\varepsilon, \eta >0$ and its propagator has $L^2$ norm equal to $1$, i.e. 
\[
||\widetilde{E}_1(t)||_{L^2 \shortrightarrow L^2}=1.
\]
Moreover, it is clear that $\widetilde{U}_1(t)=\Phi U_1(t) \Phi^{-1}$, where $U_1(t)$ is te propagator for the corresponding Wigner equation, $U_1(t) W_0^\varepsilon = \Phi^{-1} \widetilde{U}_1(t) \widetilde{W}_0^\varepsilon$.
Following  Theorem \ref{thrmWignL2}, $U_1(t)$ is an isometry in $L^2$ $\forall t \in \mathbb{R}$. (Obviously a similar statement is true for $\widetilde{U}(t):\widetilde{W}_0^\varepsilon \mapsto \widetilde{W}^\varepsilon(t)$).

Now by applying lemma \ref{lmANS} it follows that, $\forall 0\leqslant t \leqslant T$,
\begin{equation}
\begin{array}{c}
\label{eq212}
|| \rho_1^\varepsilon(t)-\widetilde{W}_1^\varepsilon(t) ||_{L^2} \leqslant
T \mathop{sup}\limits_{\tau \in [0,T]} || \widetilde{E}_1(t-\tau) F \widetilde{U}_1(\tau) \widetilde{W}^\varepsilon_0  ||_{L^2} \leqslant \\ { } \\

\leqslant T \mathop{sup}\limits_{\tau \in [0,T]} ||  F \widetilde{W}_1^\varepsilon(\tau)  ||_{L^2}.
\end{array}
\end{equation}
Recalling equation (\ref{eqF}), we are called to estimate
\begin{equation}
\label{EstF}
\begin{array}{l}
||\frac{\varepsilon\sigma_x^2}{2}\partial_x\cdot\partial_k \widetilde{W}_1^\varepsilon(\tau)+\\
+2\sum\limits_{m=2}^{\infty} 
{ 
\frac{\varepsilon^{m-1}}{(4\pi)^m} 
\sum\limits_{(m-l)mod2=1} {
 i^{l-m+1} \sigma_x^{2l} (-1)^{m-l} 
  \sum\limits_{
\begin{tiny}\begin{array}{c}
|A|=l \\
|B|=m-l
\end{array}
\end{tiny}} 
 { \frac{ \partial_x^{A+B} \widetilde{V}_1(x) }{A!B!} \partial_x^A\partial_k^B \widetilde{W}_1^\varepsilon(\tau) }  
}}||_{L^2}\leqslant \\ { } \\

\leqslant 
\frac{\varepsilon\sigma_x^2}{2}||\partial_x\cdot\partial_k \widetilde{W}_1^\varepsilon(\tau)||_{L^2}+\\
+2\sum\limits_{m=2}^{\infty} 
{ 
\frac{\varepsilon^{m-1}}{(4\pi)^m} 
\sum\limits_{(m-l)mod2=1} {
  \sigma_x^{2l}  
  \sum\limits_{
\begin{scriptsize}\begin{array}{c}
|A|=l \\
|B|=m-l
\end{array}
\end{scriptsize}} 
 { || \frac{ \partial_x^{A+B} \widetilde{V}_1(x) }{A!B!} \partial_x^A\partial_k^B \widetilde{W}_1^\varepsilon(\tau) ||_{L^2}}  
}} \leqslant 
\end{array}
\end{equation}
\[
\begin{array}{c}
\leqslant 
\frac{\varepsilon\sigma_x^2}{2}||\partial_x\cdot\partial_k \widetilde{W}_1^\varepsilon(\tau)||_{L^2}+\\
+2\sum\limits_{m=2}^{\infty} 
{ 
\frac{\varepsilon^{m-1}}{(4\pi)^m} 
\sum\limits_{(m-l)mod2=1} {
  \sigma_x^{2l}  
  \sum\limits_{
\begin{scriptsize}\begin{array}{c}
|A|=l \\
|B|=m-l
\end{array}
\end{scriptsize}} 
 { \frac{ ||\partial_x^{A+B} \widetilde{V}_1(x)||_{L^\infty} }{A!B!}
 ||  \partial_x^A\partial_k^B \widetilde{W}_1^\varepsilon(\tau) ||_{L^2}}  
}} .
\end{array}
\]

For a finite number of $m$'s, namely $m-2-\theta<-1$, assumption A1(r=0) implies that
\begin{equation}
\label{eq9io0p007}
||\partial_x^{A+B} \widetilde{V}_1(x)||_{L^\infty} \leqslant  O(1).
\end{equation}
As will be clear soon, these don't yield any intersting contribution.

Now assume that $m-2-\theta>-1$; it suffices to recall that by A1(r=0)
\begin{equation}
\label{EstVm007}
\begin{array}{c}
||\partial_x^{A+B} \widetilde{V}_1(x)||_{L^\infty} 
\leqslant  (2\pi)^m M_0 +D \frac{(2\pi)^m}{2} (\eta')^{-\frac{m-1-\theta}{2}} \Gamma\left({ \frac{m-1-\theta}{2}}\right),
\end{array}
\end{equation}
where for brevity we denote
\[
\eta'=\eta+\varepsilon\frac{\pi}{2}\sigma_x^2.
\]

To estimate $||  \partial_x^A\partial_k^B \widetilde{W}_1^\varepsilon(\tau) ||_{L^2}$ we need to remember that $|B|\geqslant 1$, therefore for each $m,l$ appearing in the sum of equation (\ref{EstF}) there is  a multi-index $|b|=1$ such that $b\leqslant B$. At this point one should recall the definition of $M_1$ (equation (\ref{eqdefm1})).  With that in mind, we have
\begin{equation}
\label{EstWderivs}
\begin{array}{c}
||  \partial_x^A\partial_k^B \widetilde{W}_1^\varepsilon(\tau) ||_{L^2}  =
(2\pi)^m ||  X^A K^B e^{-\frac{\varepsilon\pi}{2}[\sigma_x^2|X|^2+\sigma_k^2|K|^2]} \widehat{{W}}_1^\varepsilon(X,K,\tau) ||_{L^2}\leqslant \\ { } \\

\leqslant (2\pi)^m ||  X^A K^{B-b} e^{-\frac{\varepsilon\pi}{2}[\sigma_x^2|X|^2+\sigma_k^2|K|^2]}||_{L^\infty} ||K^b \widehat{{W}}_1^\varepsilon(X,K,\tau) ||_{L^2} \leqslant 
\end{array}
\end{equation}
\[
\begin{array}{c}
\leqslant (2\pi)^{(m-1) }
\left({
\prod\limits_{d=1}^n ||  X_d^{A_d} e^{-\frac{\varepsilon\pi}{2}\sigma_x^2X_d^2}||_{L^\infty} }\right)
\left({
\prod\limits_{d=1}^n ||  K_d^{(B-b)_d} e^{-\frac{\varepsilon\pi}{2}\sigma_k^2K_d^2}||_{L^\infty}  }\right)
M_1 = \\ { } \\

=(2\pi)^{(m-1) }
\prod\limits_{d=1}^n 

\left({  \frac{A_d} {  e \varepsilon \pi \sigma_x^2 } }\right)^{\frac{A_d}{2}}

\left({  \frac{(B-b)_d} {  e \varepsilon \pi \sigma_k^2 } }\right)^{\frac{(B-b)_d}{2}}
\,\,

M_1 = \\ { } \\

=M_1  \,\, (\frac{4\pi}{e\varepsilon })^{\frac{m-1}{2} } \sigma_x^{-l} \sigma_k^{-(m-l-1)}
\prod\limits_{d=1}^n A_d^{\frac{A_d}2} (B-b)_d^{\frac{(B-b)_d}2}.
\end{array}
\]
Finally, with the same arguments as above, observe that
\begin{equation}
\label{EstSK}
\begin{array}{c}
\frac{\varepsilon\sigma_x^2}{2}||\partial_x\cdot\partial_k \widetilde{W}_1^\varepsilon(\tau)||_{L^2} \leqslant 

\frac{\varepsilon\sigma_k^2}{2} 2\pi \sum\limits_{d=1}^{n}||X_d e^{-\frac{\varepsilon\pi}{2}\sigma_x^2X_d^2} ||_{L^\infty} \,\, M_1  =\\ { } \\

=\sqrt{\varepsilon} \,\, \frac{ \sqrt{\pi} \sigma_k^2 n}{\sigma_x\sqrt{e}} M_1 
\end{array}
\end{equation}

Plugging equations (\ref{EstVm}), (\ref{EstWderivs}) and (\ref{EstSK}) in (\ref{EstF}) we get 
\begin{equation}
\label{SumEstF}
\begin{array}{l}
||F \widetilde{W}_1^\varepsilon(\tau) ||_{L^2}\leqslant 

\sqrt{\varepsilon} \,\, \frac{ \sqrt{\pi} \sigma_x n}{\sqrt{e}} M_1 +  \Sigma_1 + \\ { } \\
+ \eta'^{\frac{\theta}2} \frac{M_1  D}{2} \sum\limits_{m=\lceil 1+\theta\rceil}^{\infty} 
{\left({  \frac{\varepsilon}{\eta'} \frac{\pi}{e} }\right)^{\frac{m-1}2}
\Gamma\left({ \frac{m-1-\theta}{2}}\right)} \\

{\sum\limits_{(m-l)mod2=1} {
  \sigma_x^{l}  \sigma_k^{-m+l+1}
  \sum\limits_{
\begin{scriptsize}\begin{array}{c}
|A|=l \\
|B|=m-l
\end{array}
\end{scriptsize}} 
 {  \frac{ \prod\limits_{d=1}^n A_d^{\frac{A_d}2} (B-b)_d^{\frac{(B-b)_d}2}   }{A!B!}  
}}} 
\end{array}
\end{equation}
$\Sigma_1$ stands for all the contributions from $m-2-\theta<-1$  (see equations (\ref{eq9io0p007}), (\ref{EstVm007}) ). Using our previous considerations it is easy to see that
\begin{equation}
\label{eqrrrrr00}
\begin{array}{c}
\Sigma_1= \sum\limits_{m=2}^{m_*} \frac{\varepsilon^{m-1}}{(4\pi)^m} \sum\limits_{A,B} \sigma_x^{2|A|} ||\partial_x^A\partial_k^B \widetilde{W}_1^\varepsilon||_{L^2} ||\partial_x^{A+B} \widetilde{V}_1(x) ||_{L^\infty} \leqslant \\ { } \\

\leqslant CM_1 \sum\limits_{m=2}^{m_*} \frac{\varepsilon^{\frac{m-1}{2}}}{(4\pi)^m} \sum\limits_{A,B} \sigma_x^{2|A|}  C_{A,B}= O(\sqrt{\varepsilon}M_1).
\end{array}
\end{equation}
Since we already have a term of $O(\sqrt{\varepsilon}M_1)$ -- the one coming from $\frac{\varepsilon\sigma_x^2}{2}||\partial_x\cdot\partial_k \widetilde{W}_1^\varepsilon(\tau)||_{L^2}$ -- we will not carry $\Sigma_1$ explicitly in the sequel.

Moreover, it is clear that (eventually)  $(2\pi)^m M_0 \leqslant \frac{(2\pi)^m}{2} (\eta')^{-\frac{m-1-\theta}{2}} \Gamma\left({ \frac{m-1-\theta}{2}}\right)$, so we don't treat explicitly the contributions from the first term of the rhs of (\ref{EstVm007}) for any $m$.

Now, to proceed from equation (\ref{SumEstF}) we have to understand better the inner sum. To that end, note that without loss of generality (specifically for $m$ large enough)
\begin{equation}
\label{EqGamProp}
\begin{array}{c}
A^{\frac{A}2} \approx \sqrt{ \frac{e^A}{\sqrt{2\pi A}} } \sqrt{A!}, \\ { } \\

\Gamma\left({ \frac{m-1-\theta}{2}}\right) < (2\pi)^{\frac{1}4}2^{\frac{3+2\theta}{4}} \sqrt{\frac{m!}{2^m \, m(m-1)}}.
\end{array}
\end{equation}
The first statement is the Stirling approximation; if $|A|=0$ it is clear that the contribution of the term $A^{\frac{A}2}=1$ (this can be seen very clearly if one goes back to equation (\ref{EstWderivs}) where it originates and check it. This point will always be understood as explained here, with no further explicit mention). For the proof of the second statement, see lemma \ref{LemGam}.
Now we have
\begin{equation}
\label{SumGamma}
\begin{array}{c}
\sum\limits_{(m-l)mod2=1} {
  \sigma_x^{l}  \sigma_k^{-m+l+1}
  \sum\limits_{
\begin{scriptsize}\begin{array}{c}
|A|=l \\
|B|=m-l
\end{array}
\end{scriptsize}} 
 { \frac{ \prod\limits_{d=1}^n A_d^{\frac{A_d}2} (B-b)_d^{\frac{(B-b)_d}2}   }{A!B!}  
  \Gamma\left({ \frac{m-1-\theta}{2}}\right)  }
} \approx \\ { } \\

\approx\sum\limits_{(m-l)mod2=1} {   \sigma_x^{l}  \sigma_k^{-m+l+1}
  \sum\limits_{
\begin{scriptsize}\begin{array}{c}
|A|=l \\
|B|=m-l
\end{array}
\end{scriptsize}} 
 { \frac{ \prod\limits_{d=1}^n 
 
\frac{e^{\frac{A_d+(B-b)_d}{2}} }{\sqrt{2\pi\sqrt{A_d(B-b)_d}}} \sqrt{A_d!(B-b)_d!}
   }{A!B!}  
\Gamma\left({ \frac{m-1-\theta}{2}}\right)  }
} < \\ { } \\

<\frac{(2\pi)^{\frac{1}4}2^{\frac{3+2\theta}{4}}}{(2\pi)^{\frac{n}2}} \left({\frac{e}{2}}\right)^{\frac{m-1}2}  \sum\limits_{l=0}^{m} {   \frac{\sigma_x^{l}  \sigma_k^{-m+l+1}}{\sqrt{m(m-1)}}}
  \sum\limits_{
\begin{scriptsize}\begin{array}{c}
|A|=l \\
|B|=m-l
\end{array}
\end{scriptsize}} 
 { \sqrt{\frac{ m! }{A!B!}  }  \left({ \prod_{d=1}^n \frac{1}{A_d (B-b)_d} }\right)^{\frac{1}4}
} \leqslant 
\end{array}
\end{equation}
\[
\begin{array}{c}
\leqslant
\frac{
2^{\frac{3+2\theta}{4}} \sqrt{m}} {(2\pi)^{\frac{2n-1}4} \sqrt{m(m-1)}
}
 \left({\frac{e}{2}}\right)^{\frac{m-1}2} 
\sqrt{
\sum\limits_{l=0}^{m} {

 \sigma_x^{2l}  \sigma_k^{2(l-m+1)}
  \sum\limits_{
\begin{scriptsize}\begin{array}{c}
|A|=l \\
|B|=m-l
\end{array}
\end{scriptsize}} 
 { \frac{ m! }{A!B! }  
 }} } \leqslant \\ { } \\

\leqslant\frac{2^{\frac{3+2\theta}{4}} \sigma_k }{(2\pi)^{\frac{2n-1}4}} \left({\frac{e}{2}}\right)^{\frac{m-1}2} 

\sqrt{ (n\sigma_x^2+\frac{n}{\sigma_k^2})^m  },
\end{array}
\]
where we used the observation that $ \sum\limits_{d=1}^m {x_d} \leqslant \sqrt{m} \sqrt{ \sum\limits_{d=1}^m {x^2_d}} $ and the multinomial expansion for the two last steps.

Going back to equation (\ref{SumEstF}), we get
\begin{equation}
\label{SumEstF2}
\begin{array}{c}
||F \widetilde{W}_1^\varepsilon(\tau) ||_{L^2}\leqslant 
\sqrt{\varepsilon} \,\, \frac{ \sqrt{\pi} \sigma_x n}{\sqrt{e}} M_1 +\\ { } \\
+ \eta'^{\frac{\theta}2} 

\frac{M_1  D 2^{\frac{2\theta-1}{4}}   \sqrt{n(1+\sigma_x^2\sigma_k^2)} }{(2\pi)^{\frac{2n-1}4}}
\sum\limits_{m=\lceil 1+\theta\rceil}^{\infty} 

{\left({  \frac{\varepsilon}{\eta'} \frac{n\pi (\sigma_x^2+\frac{1}{\sigma_k^2})}{2} }\right)^{\frac{m-1}2} } \leqslant \\ { } \\

\leqslant 
\sqrt{\varepsilon} \,\, \frac{ \sqrt{\pi} \sigma_x n}{\sqrt{e}} M_1 +\\ { } \\
+ \eta'^{\frac{\theta}2} 

\frac{M_1  D 2^{\frac{2\theta-1}{4}}   \sqrt{n(1+\sigma_x^2\sigma_k^2)} }{(2\pi)^{\frac{2n-1}4}}
\sum\limits_{m=2}^{\infty} 

{\left({  \frac{\varepsilon}{\eta'} \frac{n\pi (\sigma_x^2+\frac{1}{\sigma_k^2})}{2} }\right)^{\frac{m-1}2} }= 
\end{array}
\end{equation}
\[
\begin{array}{c}
=\sqrt{\varepsilon} \,\, \frac{ \sqrt{\pi} \sigma_x n}{\sqrt{e}} M_1 + \\ { } \\

+ \eta'^{\frac{\theta}2} 
\frac{M_1  D 2^{\frac{2\theta-1}{4}}   \sqrt{n(1+\sigma_x^2\sigma_k^2)} }{(2\pi)^{\frac{2n-1}4}}

\frac{\sqrt{  \frac{\varepsilon}{\eta'} \frac{n\pi (\sigma_x^2+\frac{1}{\sigma_k^2})}{2} }}{1-\sqrt{  \frac{\varepsilon}{\eta'} \frac{n\pi (\sigma_x^2+\frac{1}{\sigma_k^2})}{2} }} .
\end{array}
\]

It is now clear that the constraints on $\eta$ (recall that $\eta'=\eta+\frac{\varepsilon\pi\sigma_x^2}2$, and $\eta$ is specified in equation (\ref{SetEta})) for this to work are
\begin{equation}
\begin{array}{c}
D_1=  \frac{\varepsilon}{\eta+\frac{\varepsilon\pi\sigma_x^2}2} \frac{n\pi (\sigma_x^2+\frac{1}{\sigma_k^2})}{2}  <1 \\ { } \\

D_2=\frac{\sqrt{D_1}}{1-\sqrt{D_1}} \leqslant O(1)
\end{array}
\end{equation}
For example, by setting
\begin{equation}
\label{eqSETeta}
\eta > \varepsilon \frac{\pi}2 \cdot max\left({   (4n-1)\sigma_x^2 \,\, , \,\,   \frac{4n}{\sigma_k^2}-\sigma_x^2   }\right)
\end{equation}
it follows that
\begin{equation}
D_1 < \frac{1}{2}, \,\,\,\,\,\,\,\,\,\,\,\,\,\,\,\,\,\,\,\, D_2 <\frac{1}{\sqrt{2}-1}.
\end{equation}
In particular $\eta$ can be 
\begin{equation}
\label{SetEta}
\begin{array}{c}
\eta = \varepsilon \left({  \frac{\pi}2 \cdot max\left({   (4n-1)\sigma_x^2 \,\, , \,\,   \frac{4n}{\sigma_k^2}-\sigma_x^2   }\right) +1 }\right)=O(\varepsilon).
\end{array}
\end{equation}

Equation (\ref{eqPart1}) follows by combining equations (\ref{eq212}), (\ref{SumEstF2}) and (\ref{SetEta}).

\vskip 0.5cm
\noindent {\bf Proof of equation (\ref{eqPart2}):} It is obvious that $||\widetilde{W}_1^\varepsilon(t)-\widetilde{W}^\varepsilon(t)||_{L^2}=||\Phi({W}^\varepsilon(t)-W^\varepsilon_1(t))||_{L^2}\leqslant ||{W}^\varepsilon(t)-W_1^\varepsilon(t)||_{L^2}$.

We will use again lemma \ref{lmANS}.
Recall the definitions of $U(t)$, $U_1(t)$.  It is straightforward to compute that if
\begin{equation}
\label{eqzzvv867r}
\begin{array}{c}
LW_1^\varepsilon = \frac{2}{\varepsilon} Re \left[{ i\int{e^{2\pi i Sx}(1-e^{-\eta' S^2})\hat{V}(S)W_1^\varepsilon(x,k-\frac{\varepsilon S}{2})dS} }\right],
\end{array}
\end{equation}
it follows that
\begin{equation}
\begin{array}{c}
\mathcal{F}_{x,k \shortrightarrow X,K} \left[{ LW_1^\varepsilon }\right]=\\ { } \\
=\frac{i}{\varepsilon} \int{ \hat{W}_1^\varepsilon(X-S,K) (1-e^{-\eta' S^2})\hat{V}(S) \left[{ e^{-\pi i \varepsilon S\cdot K}-e^{\pi i \varepsilon S\cdot K}  }\right]dS  } =\\ { } \\

=2\pi  \int{ \hat{W}_1^\varepsilon(X-S,K) \,\, (1-e^{-\eta' S^2})\hat{V}(S)S\cdot K \,\, \frac{sin(\varepsilon\pi S\cdot K)}{\varepsilon\pi S\cdot K} dS  }.
\end{array}
\end{equation}
Using lemma \ref{lemyoungineq} it follows that
\begin{equation}
\begin{array}{c}
\label{eqN1}
|| LW_1^\varepsilon ||_{L^2(\mathbb{R}^{2n})}=
|| \mathcal{F} \left[{ LW_1^\varepsilon }\right] ||_{L^2(\mathbb{R}^{2n})} \leqslant \\ { } \\

\leqslant 2
\sum\limits_{i=1}^n
\pi|| K_i\hat{W}_1^\varepsilon(X,K) ||_{L^2(\mathbb{R}^{2n})} \,\,
|| (1-e^{-\eta' S^2})\hat{V}(S)S_i ||_{L^1(\mathbb{R}^{n})}  \,\,
|| \frac{sin(\varepsilon\pi S\cdot K)}{\varepsilon\pi S\cdot K} ||_{L^\infty}.
\end{array}
\end{equation}
Obviously $| \frac{sin(\varepsilon\pi S\cdot K)}{\varepsilon\pi S\cdot K} | \leqslant 1$. Moreover, recalling assumption (A3) we have
\begin{equation}
\label{eqN2}
\sum\limits_{i=1}^n|| K_i\hat{W}_1^\varepsilon(X,K) ||_{L^2(\mathbb{R}^{2n})}\leqslant M_1 \frac{1}{2\pi}.
\end{equation}
 Finally, recalling equation (\ref{eqq03}) of A1(r=0)
\begin{equation}
\begin{array}{c}
\label{eqN3}
|| (1-e^{-\eta' S^2}) S_i\hat{V}(S)  ||_{L^1(\mathbb{R}^{n})}\leqslant
  O(\eta'^{\frac{\theta}2})
\end{array}
\end{equation}
without loss of generality (i.e. $\theta<2$ without loss of generality).

Therefore summarizing equations (\ref{eqN1}), (\ref{eqN2}), (\ref{eqN3}) we finally get
\begin{equation}
\begin{array}{c}
|| LW_1^\varepsilon ||_{L^2(\mathbb{R}^{2n})}
\leqslant  M_1 C \eta'^{\frac{\theta}2}.
\end{array}
\end{equation}

Now by applying lemma \ref{lmANS} we get
\begin{equation}
\label{eq8811yyjtf}
\begin{array}{l}
W_1^\varepsilon(t)-W^\varepsilon(t)=-\int\limits_{\tau=0}^{t} {  U(t-\tau) L U_1(\tau)W^\varepsilon_0  d\tau}  \Rightarrow \\ { } \\

\Rightarrow || W_1^\varepsilon(t)-W^\varepsilon(t)||_{L^2} \leqslant\\
\,\,\,\,\,\,\,\,\,\,\,\,\,\,\,\,\,\,\,\,\,\,\,
\leqslant
 T \,\, \mathop{sup}\limits_{\tau \in [0,T]} ||U(t)||_{L^2 \shortrightarrow L^2} \,\,
\mathop{sup}\limits_{\tau \in [0,T]} ||LW_1^\varepsilon(t)||_{L^2} \leqslant \\ { } \\

\,\,\,\,\,\,\,\,\,\,\,\,\,\,\,\,\,\,\,\,\,\,\,
\leqslant
 T \,\, \left({ \eta'^{\frac{\theta}2} + \eta'}\right) \,\,C' \,\, M_1 
\end{array}
\end{equation}
for some $O(1)$ constant $C'$. 

\vskip 0.5cm
\noindent {\bf Proof of equation (\ref{eqPart3}):} Assumption A1(r=0) allows us to use
 theorem \ref{thrm1corr}, for $m=1$: there is an $O(1)$ constant $D(1,n)$ so that 
\begin{equation}
||W^\varepsilon(t)||_{H^1} \leqslant e^{t D(1,n)} ||W^\varepsilon_0||_{H^1}.
\end{equation}

Now the conclusion follows by application of lemma \ref{thrmUnSmooth} -- and more precisely its corollary described in equation (\ref{eq5rt67y}).

\vskip 0.5cm
\noindent {\bf Proof of equation (\ref{eqPart3a}):} For this we will introduce an auxiliary object. Denote $\rho_2(x,k,t)$ as the solution of
\begin{equation}
\label{2LiouvEqa123}
\begin{array}{c}
\partial_t \rho_2 +2\pi k \cdot \partial_k \rho_2 - \frac{1}{2\pi} \partial_x V \cdot \partial_k \rho_2=0,\\ { } \\
\rho_2(t=0)=\widetilde{W}^{\varepsilon}_0,
\end{array}
\end{equation}
First of all observe that, using again equation (\ref{eq5rt67y}), there is a  constant $C$, depending only on $\sigma_x,\sigma_k,n,M_0$, such that
\begin{equation}
||\rho_2(t)-\rho(t)||_{L^2}=||W_0^\varepsilon-\widetilde{W}^{\varepsilon}_0||_{L^2} \leqslant C \sqrt{\varepsilon} ||W_0^\varepsilon||_{H^1}.
\end{equation}

We will proceed in a way analogous to what we did earlier to show that $\rho$ and $\rho_1$ are close. 

By substracting equation (\ref{2LiouvEqa123}) from (\ref{SmLiouvEqa}) it follows that $h=\rho_1-\rho_2$ satisfies
\begin{equation}
\begin{array}{c}
\partial_t h +2\pi k \cdot \partial_k h - \frac{1}{2\pi} \partial_x \widetilde{V}_1 \cdot \partial_k h=\\ { } \\

=\mathcal{F}_{X,K\rightarrow x,k}^{-1} [2 \pi \int{ (1-e^{-\eta' S^2})\widehat{V}(S) \widehat{\rho}_2(X-S,K) S \cdot K dS } ], \\ { } \\

h(t=0)=0.
\end{array}
\end{equation}

By using again lemma \ref{lmANS} -- it is clear that Liouville equations with $W^{2,\infty}$ potentials are well posed in $L^2$ -- it follows that, for $t\in [0,T]$,
\begin{equation}
\begin{array}{c}
||h(t)||_{L^2} \leqslant T ||(1-e^{-\eta' S^2})|S|\widehat{V}(S)||_{L^1} M_1 \leqslant T  M_1 C(\varepsilon+\varepsilon^{\frac{\theta}2}).
\end{array}
\end{equation}

Recalling equations (\ref{eqN3}) and (\ref{eq00iibbvv}), equation (\ref{eqPart3a}) follows.

\vskip 0.5cm
The proof of theorem \ref{thrm2} is complete.

\subsection{Proof of Theorem \ref{thrm23}}
\label{sbseczoro}

The proof is to a large extent analogous to that of theorem \ref{thrm2} (using the $H^{-m}$ continuity of $U(t)$, theorem \ref{lmapqolz}, instead of the $L^2$ continuity -- theorem \ref{thrmWignL2} -- which we used before). Therefore, we will not discuss in detail points that were treated 
previously. 

The main idea is that
\[
\begin{array}{c}
||W^\varepsilon(t)-\rho(t)||_{H^{-r}} \leqslant ||  \widetilde{W}^\varepsilon(t)   -W^\varepsilon(t)||_{H^{-r}} + \\ { } \\
+||  \widetilde{W}_1^\varepsilon   -\widetilde{W}^\varepsilon||_{H^{-r}}+
||  \widetilde{W}_1^\varepsilon   -\rho_1^\varepsilon ||_{H^{-r}}+
|| \rho_1^\varepsilon(t) -\rho(t) ||_{H^{-r}}.
\end{array}
\]

\vskip 0.2cm
It will be shown that
there is a constant $C$  such that
\begin{equation}
\label{eqPart10a}
||  \widetilde{W}_1^\varepsilon   -\rho_1^\varepsilon ||_{H^{-r}} \leqslant C T || {W}_0^\varepsilon ||_{H^{-1-r}} e^{TD(r+1,n)} (\varepsilon^{\frac{\theta}2}+\sqrt{\varepsilon} ),
\end{equation}
\begin{equation}
\label{eqPart20a}
||  \widetilde{W}_1^\varepsilon   -\widetilde{W}^\varepsilon||_{H^{-r}} \leqslant CT(\varepsilon^{\frac{\theta}2}+\varepsilon) || {W}_0^\varepsilon ||_{H^{-1-r}} e^{TD(r+1,n)}  ,
\end{equation}
\begin{equation}
\label{eqPart30a}
||  \widetilde{W}^\varepsilon(t)   -W^\varepsilon(t)||_{H^{-r}} \leqslant C  || {W}_0^\varepsilon ||_{H^{1-r}} e^{TD(r-1,n)}  \sqrt{\varepsilon}, 
\end{equation}
and 
\begin{equation}
\label{eqPart40a}
\begin{array}{c}
|| \rho_1^\varepsilon(t) -\rho(t) ||_{H^{-r}} \leqslant\\ \leqslant (\varepsilon+\varepsilon^{\frac{\theta}{2}})
C T e^{T[D(r+1,n)+D(r,n)]} || W_0^\varepsilon||_{H^{-1-r}}  +\\+C \sqrt{\varepsilon} e^{T D(r,n)} || W_0^\varepsilon||_{H^{1-r}} .
\end{array}
\end{equation}

\vskip 0.5cm
Now we proceed to the proofs of the building blocks:

\noindent {\bf Proof of equation (\ref{eqPart10a}):} 

Recall that $\widetilde{E}_1(t)$ is the propagator of equation (\ref{SmLiouvEqa}), i.e. $\rho_1^\varepsilon(t)=\widetilde{E}_1(t)\widetilde{W}^\varepsilon_0$, and $\widetilde{U}_1(t)$ the propagator of equation (\ref{SWTeq1}), $\widetilde{W}_1^\varepsilon(t)=\widetilde{U}_1(t)\widetilde{W}^\varepsilon_0$. 

It follows from theorem \ref{liouv2} (applying the second sufficient condition), using the assumption A1(r), that
\[
||\widetilde{E}_1(t)||_{H^{r} \shortrightarrow H^{r}}\leqslant e^{tD(r,n)}.
\]
With respect to the sufficient condition, observe that
\[
\mathop{sup}\limits_{|a|\leqslant r+2}  ||\partial_x^aV_1||_{L^\infty} \leqslant || \widehat{V}_1(S)  |S|^{r+2}||_{L^1}  < \infty.
\]
Therefore, it follows by duality (just as in lemma \ref{lmapqolz}) that
\[
||\widetilde{E}_1(t)||_{H^{-r} \shortrightarrow H^{-r}}\leqslant e^{tD(r,n)}.
\]
(Recall that $D(r,n)$ was defined in equation (\ref{eqapp001mnj}), and can be estimated using A1(r)).

Moreover, it follows from lemma \ref{lmapqolz} that 
$U_1(t)$, $U(t)$ are bounded in $H^{-r} \shortrightarrow H^{-r}$.

Now by applying lemma \ref{lmANS} it follows that, $\forall 0\leqslant t \leqslant T$,
\begin{equation}
\begin{array}{c}
\label{eq2121}
|| \rho_1^\varepsilon(t)-\widetilde{W}_1^\varepsilon(t) ||_{H^{-r}} \leqslant
T \mathop{sup}\limits_{\tau \in [0,T]} || \widetilde{E}_1(t-\tau) F \widetilde{U}_1(\tau) \widetilde{W}^\varepsilon_0  ||_{H^{-r}} \leqslant \\ { } \\

\leqslant T e^{TD(r,n)}  \mathop{sup}\limits_{\tau \in [0,T]} ||  F \widetilde{W}_1^\varepsilon(\tau)  ||_{H^{-r}}.
\end{array}
\end{equation}
Recalling equation (\ref{eqF}), we are called to estimate
\begin{equation}
\label{EstF1}
\begin{array}{l}
||\frac{\varepsilon\sigma_x^2}{2}\partial_x\cdot\partial_k \widetilde{W}_1^\varepsilon(\tau)+\\
+2\sum\limits_{m=2}^{\infty} 
{ 
\frac{\varepsilon^{m-1}}{(4\pi)^m} 
\sum\limits_{(m-l)mod2=1} {
 i^{l-m+1} \sigma_x^{2l} (-1)^{m-l} 
  \sum\limits_{
\begin{tiny}\begin{array}{c}
|A|=l \\
|B|=m-l
\end{array}
\end{tiny}} 
 { \frac{ \partial_x^{A+B} \widetilde{V}_1(x) }{A!B!} \partial_x^A\partial_k^B \widetilde{W}_1^\varepsilon(\tau) }  
}}||_{H^{-r}}\leqslant \\ { } \\

\leqslant 
\frac{\varepsilon\sigma_x^2}{2}||\partial_x\cdot\partial_k \widetilde{W}_1^\varepsilon(\tau)||_{H^{-r}}+\\
+2\sum\limits_{m=2}^{\infty} 
{ 
\frac{\varepsilon^{m-1}}{(4\pi)^m} 
\sum\limits_{(m-l)mod2=1} {
  \sigma_x^{2l}  
  \sum\limits_{
\begin{scriptsize}\begin{array}{c}
|A|=l \\
|B|=m-l
\end{array}
\end{scriptsize}} 
 { || \frac{ \partial_x^{A+B} \widetilde{V}_1(x) }{A!B!} \partial_x^A\partial_k^B \widetilde{W}_1^\varepsilon(\tau) ||_{H^{-r}}}  
}} \leqslant \\ { } \\

\leqslant 
\frac{\varepsilon\sigma_x^2}{2}||\partial_x\cdot\partial_k \widetilde{W}_1^\varepsilon(\tau)||_{H^{-r}}+\\
+2\sum\limits_{m=2}^{\infty} 
{ 
\frac{\varepsilon^{m-1}}{(4\pi)^m} 
\sum\limits_{(m-l)mod2=1} {
  \sigma_x^{2l}  
  \sum\limits_{
\begin{scriptsize}\begin{array}{c}
|A|=l \\
|B|=m-l
\end{array}
\end{scriptsize}} 
 { \frac{ ||\partial_x^{A+B} \widetilde{V}_1(x)||_{L^\infty} }{A!B!}
 ||  \partial_x^A\partial_k^B \widetilde{W}_1^\varepsilon(\tau) ||_{H^{-r}}}  
}} .
\end{array}
\end{equation}

Now we proceed to estimate $||  \partial_x^A\partial_k^B \widetilde{W}_1^\varepsilon(\tau) ||_{H^{-r}}$:
\begin{equation}
\label{EstWderivs1}
\begin{array}{c}
||  \partial_x^A\partial_k^B \widetilde{W}_1^\varepsilon(\tau) ||_{H^{-r}}  =\\ { } \\

=(2\pi)^m ||  X^A K^B e^{-\frac{\varepsilon\pi}{2}[\sigma_x^2|X|^2+\sigma_k^2|K|^2]} \widehat{W}_1^\varepsilon(X,K,\tau) ||_{\mathcal{F} H^{-r}}\leqslant \\ { } \\

\leqslant (2\pi)^m ||  X^A K^{B-b} e^{-\frac{\varepsilon\pi}{2}[\sigma_x^2|X|^2+\sigma_k^2|K|^2]}||_{L^\infty} ||K^b \widehat{W}_1^\varepsilon(X,K,\tau) ||_{\mathcal{F}H^{-r}} \leqslant \\ { } \\

\leqslant (2\pi)^{(m-1) }
\left({
\prod\limits_{d=1}^n ||  X_d^{A_d} e^{-\frac{\varepsilon\pi}{2}\sigma_x^2X_d^2}||_{L^\infty} }\right) \\
\left({
\prod\limits_{d=1}^n ||  K_d^{(B-b)_d} e^{-\frac{\varepsilon\pi}{2}\sigma_k^2K_d^2}||_{L^\infty}  }\right)
|| \widehat{W}_1^\varepsilon(X,K,\tau) ||_{\mathcal{F}H^{-1-r}} = \\ { } \\

=(2\pi)^{(m-1) }
\prod\limits_{d=1}^n 

\left({  \frac{A_d} {  e \varepsilon \pi \sigma_x^2 } }\right)^{\frac{A_d}{2}}

\left({  \frac{(B-b)_d} {  e \varepsilon \pi \sigma_k^2 } }\right)^{\frac{(B-b)_d}{2}}
\,\,
|| \widehat{W}_0^\varepsilon ||_{\mathcal{F}H^{-1-r}} e^{tD(r+1,n)} = \\ { } \\

=|| {W}_0^\varepsilon ||_{H^{-1-r}} e^{tD(r+1,n)}  \,\, (\frac{4\pi}{e\varepsilon })^{\frac{m-1}{2} } \sigma_x^{-l} \sigma_k^{-(m-l-1)}
\prod\limits_{d=1}^n A_d^{\frac{A_d}2} (B-b)_d^{\frac{(B-b)_d}2}.
\end{array}
\end{equation}
Finally, with the same arguments as above, observe that
\begin{equation}
\label{EstSK1}
\begin{array}{c}
\frac{\varepsilon\sigma_x^2}{2}||\partial_x\cdot\partial_k \widetilde{W}_1^\varepsilon(\tau)||_{H^{-r}} \leqslant \\ { } \\

\leqslant \frac{\varepsilon\sigma_k^2}{2} 2\pi \sum\limits_{d=1}^{n}||X_d e^{-\frac{\varepsilon\pi}{2}\sigma_x^2X_d^2} ||_{L^\infty} \,\, || {W}_0^\varepsilon ||_{H^{-1-r}} e^{tD(r+1,n)}  =\\ { } \\

=\sqrt{\varepsilon} \,\, \frac{ \sqrt{\pi} \sigma_k^2 n}{\sigma_x\sqrt{e}} || {W}_0^\varepsilon ||_{H^{-1-r}} e^{tD(r+1,n)}
\end{array}
\end{equation}

Plugging equations (\ref{eqq02}), (\ref{EstWderivs1}) and (\ref{EstSK1}) in (\ref{EstF1}) we get 
\begin{equation}
\label{SumEstF1}
\begin{array}{l}
||F \widetilde{W}_1^\varepsilon(\tau) ||_{H^{-r}}\leqslant 

\sqrt{\varepsilon} \,\, \frac{ \sqrt{\pi} \sigma_x n}{\sqrt{e}}  || {W}_0^\varepsilon ||_{H^{-1-r}} e^{TD(r+1,n)}  +  \Sigma_1 + \\ { } \\
+ \eta'^{\frac{\theta}2} \frac{ || {W}_0^\varepsilon ||_{H^{-1-r}} e^{TD(r+1,n)}   D}{2} \sum\limits_{m=\lceil 1+\theta\rceil}^{\infty} 
{\left({  \frac{\varepsilon}{\eta'} \frac{\pi}{e} }\right)^{\frac{m-1}2}
\Gamma\left({ \frac{m-1-\theta}{2}}\right)} \\

{\sum\limits_{(m-l)mod2=1} {
  \sigma_x^{l}  \sigma_k^{-m+l+1}
  \sum\limits_{
\begin{scriptsize}\begin{array}{c}
|A|=l \\
|B|=m-l
\end{array}
\end{scriptsize}} 
 {  \frac{ \prod\limits_{d=1}^n A_d^{\frac{A_d}2} (B-b)_d^{\frac{(B-b)_d}2}   }{A!B!}  
}}} 
\end{array}
\end{equation}
$\Sigma_1$ stands for all the contributions from $m<\theta+1$, and behaves simlarly to $\sqrt{\varepsilon} || {W}_0^\varepsilon ||_{H^{-1-r}} e^{TD(r+1,n)}$ (see e.g. equation (\ref{eqrrrrr00}) ). 

By working in complete analogy to the proof of equation (\ref{eqPart1}), (\ref{eqPart10a}) follows.

\vskip 0.5cm
\noindent {\bf Proof of equation (\ref{eqPart20a}):} It is obvious that
\[
||\widetilde{W}_1^\varepsilon(t)-\widetilde{W}^\varepsilon(t)||_{H^{-r}} \leqslant  ||{W}_1^\varepsilon(t)-{W}^\varepsilon(t)||_{H^{-r}}.
\]
We will use again lemma \ref{lmANS} to estimate $||{W}_1^\varepsilon(t)-{W}^\varepsilon(t)||_{H^{-r}}$.
It is straightforward to compute that if
\begin{equation}
\label{eqzzvv867r1}
\begin{array}{c}
LW_1^\varepsilon = \frac{2}{\varepsilon} Re \left[{ i\int{e^{2\pi i Sx}(1-e^{-\eta' S^2})\hat{V}(S)W_1^\varepsilon(x,k-\frac{\varepsilon S}{2})dS} }\right],
\end{array}
\end{equation}
it follows that
\begin{equation}
\begin{array}{c}
\label{eqN11}
|| LW_1^\varepsilon ||_{H^{-r}}=
|| \mathcal{F} \left[{ LW_1^\varepsilon }\right] ||_{\mathcal{F}H^{-r}} \leqslant \\ { } \\

\leqslant 2
\sum\limits_{i=1}^n
\pi|| K_i\hat{W}_1^\varepsilon(X,K,t) ||_{\mathcal{F}H^{-r}} \,\,
|| (1-e^{-\eta' S^2})\hat{V}(S)S_i ||_{L^1}  \,\,
|| \frac{sin(\varepsilon\pi S\cdot K)}{\varepsilon\pi S\cdot K} ||_{L^\infty} \leqslant \\ { } \\

\leqslant 2
n
\pi  || {W}_0^\varepsilon ||_{H^{-1-r}} e^{TD(r+1,n)}  \,\,
M_3.
\end{array}
\end{equation}

Recalling equation (\ref{eqq03}) of A1(r=0), $|| (1-e^{-\eta' S^2})\hat{V}(S)S_i ||_{L^1}  =O(\varepsilon + \varepsilon^{\frac{\theta}{2}})$;
equation (\ref{eqPart20a}) now follows. 

\vskip 0.5cm
\noindent {\bf Proof of equation (\ref{eqPart30a}):} It follows readily by application of lemma \ref{thrmUnSmooth} for $p=2$, $m=-r$, $s=2$.

\vskip 0.5cm
\noindent {\bf Proof of equation (\ref{eqPart40a}):} Recall the definition of $\rho_2(x,k,t)$, equation (\ref{2LiouvEqa123}).

First of all observe that, using again lemma \ref{thrmUnSmooth}, there is a  constant $C$, depending only on $\sigma_x,\sigma_k,n,M_0$, such that
\begin{equation}
||\rho_2(t)-\rho(t)||_{H^{-r}}\leqslant e^{T D(r,n)} ||W_0^\varepsilon-\widetilde{W}^{\varepsilon}_0||_{H^{-r}} \leqslant C \sqrt{\varepsilon} e^{T D(r,n)} || W_0^\varepsilon||_{H^{1-r}}.
\end{equation}

By substracting equation (\ref{2LiouvEqa123}) from (\ref{SmLiouvEqa}) it follows that $h=\rho_1-\rho_2$ satisfies
\begin{equation}
\label{eqsustraction}
\begin{array}{c}
\partial_t h +2\pi k \cdot \partial_k h - \frac{1}{2\pi} \partial_x \widetilde{V}_1 \cdot \partial_k h=\\ { } \\

=\mathcal{F}_{X,K\rightarrow x,k}^{-1} [2\pi \int{ (1-e^{-\eta' S^2})\widehat{V}(S) \widehat{\rho}_2(X-S,K) S \cdot K dS } ], \\ { } \\

h(t=0)=0.
\end{array}
\end{equation}

By using again lemmata \ref{lmANS}, \ref{lemyoungineq}  it follows that, for $t\in [0,T]$,
\begin{equation}
\begin{array}{c}
||h(t)||_{H^{-r}} 
\leqslant C T e^{TD(r,n)} ||(1-e^{-\eta' S^2})|S|\widehat{V}(S)||_{L^1} ||\rho_2(t)||_{H^{-1-r}} \leqslant \\ { } \\
\leqslant
(\varepsilon+\varepsilon^{\frac{\theta}{2}}) C' T e^{T[D(r+1,n)+D(r,n)]} || W_0^\varepsilon||_{H^{-1-r}}   .
\end{array}
\end{equation}
for some constant $C'$ depending only on $n,\sigma_x,\sigma_k,M_0$. (We also used theorem \ref{liouv2} to estimate the rhs of $||h(t)||_{H^{-r}} \leqslant ||\rho_1(t)||_{H^{-r}}+||\rho_2(t)||_{H^{-r}}$).

Using the triangle inequality, equation (\ref{eqPart40a}) follows.

\vskip 0.5cm
The proof of theorem \ref{thrm23} is complete.

\subsection{Proof of Theorem \ref{thrm24}}

The proof is very similar to that of theorem \ref{thrm23}. We will only discuss explicitly the points where differences arise.

Of course the main idea is that
\[
\begin{array}{c}
||W^\varepsilon(t)-\rho(t)||_{H^{r}} \leqslant ||  \widetilde{W}^\varepsilon(t)   -W^\varepsilon(t)||_{H^{r}} + \\ { } \\
+||  \widetilde{W}_1^\varepsilon   -\widetilde{W}^\varepsilon||_{H^{r}}+
||  \widetilde{W}_1^\varepsilon   -\rho_1^\varepsilon ||_{H^{r}}+
|| \rho_1^\varepsilon(t) -\rho(t) ||_{H^{r}}.
\end{array}
\]
It will be shown that
there is a constant $C$  such that
\begin{equation}
\label{eqPart10abb}
||  \widetilde{W}_1^\varepsilon   -\rho_1^\varepsilon ||_{H^{r}} \leqslant C T|| {W}_0^\varepsilon ||_{H^{1+r}} e^{TD(r+1,n)} (\varepsilon^{\frac{\theta}2}+\sqrt{\varepsilon} ),
\end{equation}
\begin{equation}
\label{eqPart20ab}
||  \widetilde{W}_1^\varepsilon   -\widetilde{W}^\varepsilon||_{H^{r}} \leqslant CT  (\varepsilon+\varepsilon^\frac{\theta-r}{2}) || {W}_0^\varepsilon ||_{H^{1+r}} e^{T[D(r+1,n)+D(r,n)]}  ,
\end{equation}
\begin{equation}
\label{eqPart30ab}
||  \widetilde{W}^\varepsilon(t)   -W^\varepsilon(t)||_{H^{r}} \leqslant C  || {W}_0^\varepsilon ||_{H^{1+r}} e^{TD(r+1,n)}  \sqrt{\varepsilon}, 
\end{equation}
and 
\begin{equation}
\label{eqPart40ab}
\begin{array}{c}
|| \rho_1^\varepsilon(t) -\rho(t) ||_{H^{r}} 
\leqslant C  \sqrt{\varepsilon} e^{T D(r,n)} || W_0^\varepsilon||_{H^{1+r}} + \\
+C T e^{T[D(r,n)+D(r+1,n)]} || W_0^\varepsilon||_{H^{r+1}} (\varepsilon+ \varepsilon^\frac{\theta-r}{2}).
\end{array}
\end{equation}

\vskip 0.5cm
Now we proceed to the proofs of the building blocks:

\noindent {\bf Proof of equation (\ref{eqPart10abb}):} It is entirely analogous with the proof of equation (\ref{eqPart10a}), in the proof of theorem \ref{thrm23} (section \ref{sbseczoro}) -- it is not necessary to replicate it here.

\vskip 0.5cm
\noindent {\bf Proof of equation (\ref{eqPart20ab}):} It is obvious that
\[
||\widetilde{W}_1^\varepsilon(t)-\widetilde{W}^\varepsilon(t)||_{H^{r}} \leqslant  ||{W}_1^\varepsilon(t)-{W}^\varepsilon(t)||_{H^{r}}.
\]
We will use again lemma \ref{lmANS} to estimate $||{W}_1^\varepsilon(t)-{W}^\varepsilon(t)||_{H^{r}}$.
Using lemma \ref{lemyoungineq} it follows that, if $L$ is defined as in equation (\ref{eqzzvv867r1}),
then
\begin{equation}
\begin{array}{c}
\label{eqN11b}
|| LW_1^\varepsilon ||_{H^{r}}=
|| \mathcal{F} \left[{ LW_1^\varepsilon }\right] ||_{\mathcal{F}H^{r}} \leqslant \\ { } \\

\leqslant C
|| \,\,|K|\hat{W}_1^\varepsilon(X,K,t) ||_{\mathcal{F}H^{r}} \,\,
|| (1-e^{-\eta' S^2}) (|S|^{r+1}+|S|) \hat{V}(S) ||_{L^1}  \,\,
|| \frac{sin(\varepsilon\pi S\cdot K)}{\varepsilon\pi S\cdot K} ||_{L^\infty} \leqslant \\ { } \\

\leqslant C
  || {W}_0^\varepsilon ||_{H^{1+r}} e^{TD(r+1,n)}  \,\,
|| (1-e^{-\eta' S^2})(|S|^{r+1}+|S|)\hat{V}(S)  ||_{L^1}  .
\end{array}
\end{equation}
Now recall assumption A1(r) to see that 
\begin{equation}
\begin{array}{c}
|| (1-e^{-\eta' S^2}) |S|^{r+1}\hat{V}(S)  ||_{L^1(\mathbb{R}^{n})}\leqslant C (\varepsilon+\varepsilon^{\frac{\theta-r}{2}})
\end{array}
\end{equation}
Therefore 
\begin{equation}
\begin{array}{c}
|| LW_1^\varepsilon ||_{H^{r}}
\leqslant  C (\varepsilon+\varepsilon^\frac{\theta-r}{2}) || {W}_0^\varepsilon ||_{H^{1+r}} e^{TD(r+1,n)}  
\end{array}
\end{equation}

Equation (\ref{eqPart20ab}) follows by applying lemma \ref{lmANS}. 

\vskip 0.5cm
\noindent {\bf Proof of equation (\ref{eqPart30ab}):} It follows readily by application of lemma \ref{thrmUnSmooth} for $p=2$, $m=-r$, $s=2$.

\vskip 0.5cm
\noindent {\bf Proof of equation (\ref{eqPart40ab}):} Recall the definition of $\rho_2(x,k,t)$, equation (\ref{2LiouvEqa123}).

First of all observe that, using again lemma \ref{thrmUnSmooth}, there is a  constant $C$, depending only on $\sigma_x,\sigma_kn,M_0$, such that
\begin{equation}
||\rho_2(t)-\rho(t)||_{H^{r}}\leqslant e^{T D(r,n)} ||W_0^\varepsilon-\widetilde{W}^{\varepsilon}_0||_{H^{r}} \leqslant C \sqrt{\varepsilon} e^{T D(r,n)} || W_0^\varepsilon||_{H^{1+r}}.
\end{equation}

Applying lemma \ref{lemyoungineq} to equation (\ref{eqsustraction}) once more (and using again A1(r) ), it follows that for $t\in [0,T]$
\begin{equation}
\begin{array}{c}
||h(t)||_{H^{r}} \leqslant C   ||(1-e^{-\eta' S^2})(|S|+|S|^{r+1})\widehat{V}(S)||_{L^1} ||\rho_2(t)||_{H^{r+1}} \leqslant \\ { } \\
\leqslant 
  C  e^{TD(r+1,n)} || W_0^\varepsilon||_{H^{r+1}}  (\varepsilon+\varepsilon^{\frac{\theta-r}{2}}).
\end{array}
\end{equation}

Using the triangle inequality and lemma \ref{lmANS}, equation (\ref{eqPart40ab}) follows.

\vskip 0.5cm
The proof of theorem \ref{thrm24} is complete.

\subsection{Auxiliary lemmata} \label{ecccqq}

\noindent {\bf Proof of lemma \ref{lemmqq}}: 

First of all
\[
\begin{array}{c}
\int{ |\widehat{V}(S)|\,|S|^{r+2}dS }\leqslant \int\limits_{|S|<1}{ |\widehat{V}(S)|dS } + \int\limits_{|S|>1}{ |\widehat{V}(S)|\,|S|^{r+2}dS } \leqslant \\ { } \\

\leqslant  M_0+ D\int\limits_{s=1}^{\infty} { s^{r+2+n-1-n-1-\theta}ds } < \infty
\end{array}
\]
because $r-\theta<-1$. This, together with equation (\ref{eqmklop007}), implies equation (\ref{eqmklop00}).

\vskip 0.15cm
Now denote for convenience
\begin{equation}
\label{eqEtaP}
\eta'=\eta+\frac{\varepsilon\pi \sigma_x^2}2.
\end{equation}

One observes that
\begin{equation}
\label{EstVm0}
\begin{array}{c}
||\partial_x^{A+B} \widetilde{V}_1(x)||_{L^\infty} \leqslant (2\pi)^m \int{ |\hat{V}(S) \left({ \prod\limits_{d=1}^nS_d^{(A+B)_d} }\right) e^{-\eta' |S|^2}|dS }
\leqslant  \\ 

\leqslant (2\pi)^m \int{ |\hat{V}(S)| |S|^m e^{-\eta' |S|^2}dS }
\leqslant  \\ 

\leqslant (2\pi)^m  \left({  \int\limits_{|S|\leqslant 1}{ |\hat{V}(S)| |S|^m dS }+ D \int\limits_{r=1}^\infty{ r^{(m-2)-\theta} e^{-\eta' r^2}dr } }\right).
\end{array}
\end{equation}
For a finite number of $m$'s, $m-2-\theta<-1$, in which case, using equation (\ref{EstVm0}) and equation (\ref{eqmklop00}), for $|A|+|B|\leqslant m$ we have
\begin{equation}
\label{eq9io0p}
||\partial_x^{A+B} \widetilde{V}_1(x)||_{L^\infty} \leqslant  (2\pi)^m \left({ M_0 + D \int\limits_{r=1}^\infty{ r^{(m-2)-\theta} dr } }\right) =O(1)
\end{equation}
 Now assume that $m-2-\theta>-1$;\footnote{We don't have to worry about $m-2-\theta=-1$ because we have restricted $\theta \notin \mathbb{N}$.}
then equation (\ref{EstVm0}) leads to
\begin{equation}
\label{EstVm}
\begin{array}{c}
||\partial_x^{A+B} \widetilde{V}_1(x)||_{L^\infty} 
<  (2\pi)^m M_0 +D (2\pi)^m\int\limits_{r=0}^\infty{ r^{(m-2)-\theta} e^{-\eta' S^2}dS }= \\

=  (2\pi)^m M_0 +D \frac{(2\pi)^m}{2} (\eta')^{-\frac{m-1-\theta}{2}} \Gamma\left({ \frac{m-1-\theta}{2}}\right).
\end{array}
\end{equation}

This concludes the proof of equations (\ref{eqq01}) and (\ref{eqq02}).

\vskip 0.5cm
Now observe that
\begin{equation}
\begin{array}{c}
\label{eqN3b}
|| (1-e^{-\eta' S^2}) |S|^{r+1}\hat{V}(S)  ||_{L^1(\mathbb{R}^{n})}\leqslant\\ { } \\

\leqslant\int\limits_{|S|\leqslant 1} {(1-e^{-\eta' S^2}) |S|^{r+1} |\hat{V}(S)| dS} + 
D \int\limits_{|S|>1}{|S|^{r+1} \frac{1-e^{-\eta' |S|^2}}{|S|^{n+1+\theta}} dS}\leqslant \\ { } \\

\leqslant \eta' M_0 + 

D \int\limits_{s=1}^{\infty}{s^{r+1} \frac{1-e^{-\eta' s^2}}{s^{n+1+\theta}} s^{n-1}ds} \leqslant \\ { } \\

\leqslant \eta' M_0+  D \left[{ \int\limits_{1}^{A} { s^{r-1-\theta}\eta' s^2 ds}+
\int\limits_{s=A}^{\infty} { s^{r-1-\theta}  ds}
 }\right] =\\ { } \\

=\eta' M_0 + D\left({ \eta' \frac{A^{r+2-\theta}-1}{r+2-\theta} + \frac{A^{r-\theta}}{\theta-r}  }\right)=:M_3.
\end{array}
\end{equation}
By selecting $A=\eta'^{-\frac{1}{2}}$, we get
\begin{equation}
\label{eqN4b}
\begin{array}{c}
M_3(r)=\eta' M_0 + D\left({ \eta' \frac{(\eta'^{-\frac{1}{2}})^{r+2-\theta}-1}{r+2-\theta} + \frac{(\eta'^{-\frac{1}{2}})^{r-\theta}}{\theta-r}  }\right)=\\ { } \\

=\eta' M_0
+D\frac{ \eta'^{  \frac{\theta-r}{2}  }-\eta'  }
{r+2-\theta}

+\frac{\eta'^{\frac{\theta-r}{2}}}{\theta-r} .
\end{array}
\end{equation}

By recalling our assumption that $\theta > r+1$ it is easy to see that
\begin{equation}
\label{eqN412b}
M_3(r)\leqslant C (\varepsilon^\frac{\theta-r}{2}+\varepsilon)=o(1).
\end{equation}

\vskip 0.5cm
This concludes the proof of lemma \ref{lemmqq}.

\begin{lemma}[A version of the Young inequality] \label{lemyoungineq} Let $f\in L^1(\mathbb{R}^n)$. Then,
if $g \in L^2(\mathbb{R}^{2n})$,
\begin{equation}
\label{eqtraul1}
|| \int{f(s)g(x-s,k)ds} ||_{L^2(\mathbb{R}^{2n})} \leqslant ||f||_{L^1(\mathbb{R}^{n})} ||g||_{L^2(\mathbb{R}^{2n})} ;
\end{equation}
if $\widehat{g} \in H^{-m}(\mathbb{R}^{2n})$,
\begin{equation}
\label{eqtraul2}
|| \int{f(s)g(x-s,k)ds} ||_{\mathcal{F} H^{-m}(\mathbb{R}^{2n})} \leqslant ||f||_{L^1(\mathbb{R}^{n})} ||g||_{\mathcal{F}H^{-m}(\mathbb{R}^{2n})} ;
\end{equation}
and if $\widehat{g} \in H^{m}(\mathbb{R}^{2n})$, $f(s)(1+|s|^m)\in L^1(\mathbb{R}^n)$,
\begin{equation}
\label{eqtraul3}
|| \int{f(s)g(x-s,k)ds} ||_{\mathcal{F}H^m(\mathbb{R}^{2n})} \leqslant C ||g||_{\mathcal{F}H^m(\mathbb{R}^{2n})} \sum\limits_{l=0}^m ||f(s)|s|^l||_{L^1(\mathbb{R}^{n})} .
\end{equation}
\end{lemma}

\noindent {\bf Proof of equation (\ref{eqtraul1}):}
\[
\begin{array}{c}
|| \int{f(s)g(x-s,k)ds} ||_{L^2(\mathbb{R}^{2n})}^2 = \int{ f(s)g(x-s,k) \bar{f}(y) \bar{g}(x-y,k) dxdk}=\\
= \int{ \int{g(x-s,k)  \bar{g}(x-y,k)dxdk} f(s) \bar{f}(y)dsdy}\leqslant ||f||_{L^1(\mathbb{R}^{n})}^2 ||g||_{L^2(\mathbb{R}^{2n})}^2.
\end{array}
\]

\noindent {\bf Proof of equation (\ref{eqtraul2}):} Denote $\mathcal{B}=\{ \phi | \,\,||\phi||_{\mathcal{F}H^m}\leqslant 1 \}$. Then
\[
\begin{array}{c}
|| \int{f(s)g(x-s,k)ds} ||_{H^{-m}(\mathbb{R}^{2n})} = \mathop{sup}\limits_{\phi \in \mathcal{B}} \int{ f(s)g(x-s,k)ds\phi(x,k)dxdk } \leqslant \\
\leqslant ||f||_{L^1(\mathbb{R}^{n})} \mathop{sup}\limits_{\phi \in \mathcal{B},s} \int{ g(x-s,k)\phi(x,k)dxdk } \leqslant ||f||_{L^1(\mathbb{R}^{n})} ||g||_{H^{-m}(\mathbb{R}^{2n})}.
\end{array}
\]

\noindent {\bf Proof of equation (\ref{eqtraul3}):} Recall that
\[
\begin{array}{c}
|| \int{f(s)g(x-s,k)ds} ||_{\mathcal{F}H^m(\mathbb{R}^{2n})} =\\
=\sum\limits_{|A|+|B|\leqslant m} || (2\pi x)^A (2\pi k)^B \int{f(s)g(x-s,k)ds} ||_{L^2(\mathbb{R}^{2n})},
\end{array}
\]
therefore it suffices to work for one choice of $A,B$ with $|A|+|B|\leqslant m$:
\[
\begin{array}{c}
|| \int{f(s)g(x-s,k)ds}x^A k^B ||_{L^2(\mathbb{R}^{2n})} = \\
=|| \int{f(s)g(x-s,k)ds}(s+(x-s))^A k^B ||_{L^2(\mathbb{R}^{2n})} \leqslant \\
\leqslant \sum\limits_{l=0}^A \binom{A}{l} ||f(s)|s|^{A-l}||_{L^1(\mathbb{R}^{n})} ||g(x,k)x^lk^B||_{L^2(\mathbb{R}^{2n})} \leqslant \\ \leqslant
||g||_{H^m(\mathbb{R}^{2n})}  \sum\limits_{l=0}^A \binom{A}{l} ||f(s)|s|^{|A-l|}||_{L^1(\mathbb{R}^{n})} .
\end{array}
\]

\begin{lemma}[An observation on perturbations] \label{lmANS} Consider two equations,
\begin{equation}
\label{eqa1}
\begin{array}{c}
v_t+Tv=0,\\
v(t=0)=w_0,
\end{array}
\end{equation}
and
\begin{equation}
\label{eqaa1}
\begin{array}{c}
u_t+Tu=Lu,\\
u(t=0)=w_0.
\end{array}
\end{equation}
Denote $E$ and $U$ their propagators respectively, i.e. $v(t)=E(t)w_0$, $u(t)=U(t)w_0$. Then the difference of the solutions can be cast as
\begin{equation}
\label{eqansa}
u(t)-v(t)=-\int\limits_{\tau=0}^t {  E(t-\tau)LU(\tau)w_0  d\tau}=\int\limits_{\tau=0}^t {  U(t-\tau)LE(\tau)w_0  d\tau}
\end{equation}

Moreover, if
\begin{equation}
\label{eqaa12}
\begin{array}{c}
r_t+Tr=f(t),\\
r(t=0)=w_0,
\end{array}
\end{equation}
then
\begin{equation}
\label{eqansb}
r(t)-v(t)=-\int\limits_{\tau=0}^t{E(t-\tau)f(\tau)d\tau}.
\end{equation}

\end{lemma}

\noindent {\bf Remark:} The considerations here are formal; questions of well-posedness of the equations etc should be considered each time the lemma is applied.

\noindent {\bf Proof of equation (\ref{eqansa}):} Set
\begin{equation}
\label{eq35}
\psi(t)=U(-t)(u(t)-v(t)).
\end{equation}
Then
\begin{equation}
\begin{array}{c}
U(t)\psi(t)=u(t)-v(t) \,\, \Rightarrow \,\, \dot{U}(t)\psi(t)+U(t)\dot{\psi}(t)=\dot{u}(t)-\dot{v}(t) \,\, \Rightarrow \\

\Rightarrow \,\, -(T+L)U(t)\psi(t)+U(t)\dot{\psi}(t)=-Tu(t)-Lu(t)+Tv(t) \,\, \Rightarrow \\
\Rightarrow \,\, -(T+L)(u(t)-v(t))+U(t)\dot{\psi}(t)=-Tu(t)-Lu(t)+Tv(t) \,\, \Rightarrow \\

\Rightarrow \,\, \dot{\psi}(t)=-U(-t)Lv(t) \,\, \Rightarrow \,\, U(t)\psi(t)=-U(t)\int\limits_{\tau=0}^t { U(-\tau)LE(\tau)w_0 d\tau} \,\, \Rightarrow\\

\Rightarrow u(t)-v(t)=-\int\limits_{\tau=0}^t {U(t-\tau)LE(\tau)w_0  d\tau}.
\end{array}
\end{equation}

The other equality of equation (\ref{eqansa}) follows essentially in the same way. The proof is complete.

\noindent {\bf Proof of equation (\ref{eqansb}):} Set $r'(t)=v(t)-\int\limits_{\tau=0}^t{E(t-\tau)f(\tau)d\tau}$; by straightforward computation it follows that $\partial_t r' +Tr'=f(t)$.

The proof is complete.

\begin{lemma}[Interchanging summation and integration]
\label{lemJusInter}
We will show that
\[
\begin{array}{c}
\int{e^{2\pi i Sx-\frac{\varepsilon\pi}{2}\sigma_x^2S^2}\hat{V}_1(S)\widetilde{W}_1^\varepsilon(x+\frac{i\varepsilon\sigma_x^2S}{2},k-\frac{\varepsilon S}{2})dS}= \\ { } \\

=\int{e^{2\pi i Sx-(\eta+\frac{\varepsilon\pi}{2}\sigma_x^2)S^2}\hat{V}(S)
\sum\limits_{m=0}^{+\infty} \frac{\varepsilon^m}{m!2^m} [ \sum\limits_{d=1}^n S_d(i\sigma_x^2\partial_{x_d}-\partial_{k_d})]^m
\widetilde{W}_1^\varepsilon(x,k)dS}= \\ { } \\

=\sum\limits_{m=0}^{+\infty}
\int{e^{2\pi i Sx-(\eta+\frac{\varepsilon\pi}{2}\sigma_x^2)S^2}\hat{V}(S)
\frac{\varepsilon^m}{m!2^m} [ \sum\limits_{d=1}^n S_d(i\sigma_x^2\partial_{x_d}-\partial_{k_d})]^m
\widetilde{W}_1^\varepsilon(x,k)dS},
\end{array}
\]
i.e. that interchanging summation and integration is justified.
\end{lemma}

\noindent {\bf Proof:} Denote
\[
g(S)=e^{-\frac{\varepsilon\pi}{2}\sigma_x^2S^2}\widetilde{W}_1^\varepsilon(x+\frac{i\varepsilon\sigma_x^2S}{2},k-\frac{\varepsilon S}{2}),
\]
\[
g_m(S)=e^{-(\eta+\frac{\varepsilon\pi}{2}\sigma_x^2)S^2}  \sum\limits_{m=0}^{+\infty} \frac{\varepsilon^m}{m!2^m} [ \sum\limits_{d=1}^n S_d(i\sigma_x^2\partial_{x_d}-\partial_{k_d})]^m
\widetilde{W}_1^\varepsilon(x,k).
\]

We know that $\sum\limits_m g_m(S)=g(S)$ pointwise (even uniformly on compacts), since it is the Taylor expansion of an entire analytic function. 

Observe that since the potential always satisfies A1(r=0),  
\begin{equation}
\widehat{V}(S) \in L^1;
\end{equation}
therefore, it suffices to show that
\[
\sum\limits_{m} |g_m(S)| \leqslant C,
\]
the result then follows by dominated convergence. Indeed,
\begin{equation}
\label{eqtsattsa}
\begin{array}{c}
|g_m(S)| \leqslant\\ { } \\

\leqslant e^{-(\eta+\frac{\varepsilon\pi}{2}\sigma_x^2)S^2}  \sum\limits_{m=0}^{+\infty} \frac{\varepsilon^m}{m!2^m} 
 \sum\limits_{|A|=m} \frac{m!}{A_1!...A_n!} \left|{ \prod\limits_{d=1}^n \left[{ S_d(i\sigma_x^2\partial_{x_d}-\partial_{k_d})}\right]^{A_d}
\widetilde{W}_1^\varepsilon(x,k) }\right| \leqslant \\ { } \\

\leqslant e^{-(\eta+\frac{\varepsilon\pi}{2}\sigma_x^2)S^2}  \sum\limits_{m=0}^{+\infty} \frac{\varepsilon^m}{m!2^m} 
 \sum\limits_{|A|=m} \binom{m}{A} S^A 
 \prod\limits_{d=1}^n \sum\limits_{l=0}^{A_d} \binom{A_d}{l} \left|{(i\sigma_x^2\partial_{x_d})^{l}(-\partial_{k_d})^{A_d-l}
\widetilde{W}_1^\varepsilon(x,k) }\right|.
\end{array}
\end{equation}

Here observe that
\[
\begin{array}{c}
\left|{(i\sigma_x^2\partial_{x_d})^{l}(-\partial_{k_d})^{A_d-l}
\widetilde{W}_1^\varepsilon(x,k) }\right| \leqslant \\ { } \\

\leqslant  max \{1,\sigma_x^{2l} \}  (2\pi)^m \int{ X_d^lK_d^{A_d-l} e^{-\frac{\varepsilon\pi}{2}(\sigma_x^2|X|^2+\sigma_k^2|K|^2)} |\widehat{W}_1(X,K)| dXdK} \leqslant \\ { } \\

\leqslant  max \{1,\sigma_x^{2l} \}  (2\pi)^m || X_d^l e^{-\frac{\varepsilon\pi}{2}\sigma_x^2|X|^2} ||_{L^2} || K_d^{A_d-l} e^{-\frac{\varepsilon\pi}{2}\sigma_k^2|K|^2} ||_{L^2}   || \widehat{W}_1 ||_{L^2}.
\end{array}
\]
Moreover,
\[
\begin{array}{c}
|| X_d^l e^{-\frac{\varepsilon\pi}{2}\sigma_x^2|X|^2} ||_{L^2}^2 = (\int{ e^{-\varepsilon\pi\sigma_x^2x^2}dx })^{n-1} 
\int{ x^{2l} e^{-\varepsilon\pi\sigma_x^2 x^2}dx } = \\ { } \\

= ( \varepsilon \sigma_x^2)^\frac{1-n}{2}  \frac{ \Gamma \left({ \frac{2l-1}{2}  }\right)   }{(\pi\varepsilon\sigma_x^2)^{\frac{2l-1}{2}}},
\end{array}
\]
and similarly
\[
\begin{array}{c}
|| K_d^{A_d-l} e^{-\frac{\varepsilon\pi}{2}\sigma_k^2|K|^2} ||_{L^2}^2 
= ( \varepsilon \sigma_k^2)^\frac{1-n}{2}  \frac{ \Gamma \left({ \frac{2(A_d-l)-1}{2}  }\right)   }{(\pi\varepsilon\sigma_k^2)^{\frac{2(A_d-l)-1}{2}}},
\end{array}
\]
and therefore
\[
\begin{array}{c}
\left|{(i\sigma_x^2\partial_{x_d})^{l}(-\partial_{k_d})^{A_d-l}
\widetilde{W}_1^\varepsilon(x,k) }\right| \leqslant \\ { } \\

\leqslant  max \{1,\sigma_x^{2l} \}  (2\pi)^m ( \varepsilon \sigma_x \sigma_k)^\frac{1-n}{2} 
(\varepsilon\pi)^{\frac{1-m}{2}} \sigma_x^{\frac{1-2l}2} \sigma_k^{\frac{1-2(A_d-l)}2} \\

\sqrt{  \Gamma \left({ \frac{2l-1}{2}  }\right) \Gamma \left({ \frac{2(A_d-l)-1}{2}  }\right)}

|| \widehat{W}_1 ||_{L^2} \leqslant \\ { } \\

\leqslant C  max \{1,\sigma_x^{2l} \}
max \{ 1, \sigma_x^{\frac{1-2l}2} \} max\{ 1, \sigma_k^{\frac{1-2(A_d-l)}2} \} \\

  (2\pi)^m 
(\varepsilon\pi)^{\frac{1-m}{2}} 

\sqrt{  \Gamma \left({ \frac{2l-1}{2}  }\right) \Gamma \left({ \frac{2(A_d-l)-1}{2}  }\right)} \leqslant \\ { }\\

\leqslant C  B^m
\varepsilon^{\frac{1-m}{2}} 

\sqrt{  \Gamma \left({ \frac{2l-1}{2}  }\right) \Gamma \left({ \frac{2(A_d-l)-1}{2}  }\right)}
\end{array}
\]
where $C,B<\infty$ are independent of $m$. Plugging this back into equation (\ref{eqtsattsa}), there are (possibly new) constants $C,B<\infty$ such that
\begin{equation}
\label{eqtsattsa1}
\begin{array}{c}
|g_m(S)| 

\leqslant C e^{-(\eta+\frac{\varepsilon\pi}{2}\sigma_x^2)S^2}  \sum\limits_{m=0}^{+\infty} B^m \frac{\varepsilon^m}{m!} \\

\,\,\,\,\,\,\,\,\,\,\,\,\,\,\,\,\,\,\,\,\,\,\,\,\,\,\,\,\,\,\,\,\,\,\,\,
 \sum\limits_{|A|=m} \binom{m}{A} S^A 
 \prod\limits_{d=1}^n \sum\limits_{l=0}^{A_d} \binom{A_d}{l} \varepsilon^{\frac{1-m}{2}} 
\sqrt{  \Gamma \left({ \frac{2l-1}{2}  }\right) \Gamma \left({ \frac{2(A_d-l)-1}{2}  }\right)} \leqslant \\ { } \\

\leqslant C' e^{-(\eta+\frac{\varepsilon\pi}{2}\sigma_x^2)S^2}  \sum\limits_{m=0}^{+\infty} B^m \varepsilon^{\frac{m}{2}}
 \sum\limits_{|A|=m}  S^A 
 \prod\limits_{d=1}^n \sum\limits_{l_d=0}^{A_d}  
\frac{\sqrt{  \Gamma \left({ \frac{2l_d-1}{2}  }\right)  \Gamma \left({ \frac{2(A_d-l_d)-1}{2}  }\right)}}{l_d!(A_d-l_d)!} < \\ { } \\

<C' e^{-(\eta+\frac{\varepsilon\pi}{2}\sigma_x^2)S^2}  \sum\limits_{m=0}^{+\infty} B^m \varepsilon^{\frac{m}{2}}
 \sum\limits_{|A|=m}  S^A 
 \prod\limits_{d=1}^n \sum\limits_{l_d=0}^{A_d}  
\frac{\sqrt{ (l_d-1)!(A_d-l_d-1)!}}{l_d!(A_d-l_d)!} .
\end{array}
\end{equation}

At this point, using the Cauchy-Schwarz inequality one observes that
\[
\begin{array}{c}
 \prod\limits_{d=1}^n \sum\limits_{l_d=0}^{A_d}  
\frac{\sqrt{ (l_d-1)!(A_d-l_d-1)!}}{l_d!(A_d-l_d)!} =

 \prod\limits_{d=1}^n \sqrt{ \sum\limits_{l_d=0}^{A_d}  
 \frac{ A_d!}{l_d!(A_d-l_d)!} } \frac{1}{\sqrt{A_d!}} 
\sqrt{ \sum\limits_{l_d=1}^{A_d-1}  \frac{1}{l_d(A_d-l_d)} }
  \leqslant \\ { } \\

\leqslant  \prod\limits_{d=1}^n \sqrt{\frac{A_d}{A_d!(A_d-1)}} \sqrt{ \sum\limits_{l_d=0}^{A_d}  
 \frac{ A_d!}{l_d!(A_d-l_d)!} } \approx

\sqrt{\frac{2^{A_d}}{A_d!}}  .
\end{array}
\]
Finally, observe that
\[
\left|{  e^{-(\eta+\frac{\varepsilon\pi}{2}\sigma_x^2)S^2} S^A }\right| \leqslant B^m (\eta+\frac{\varepsilon\pi}{2}\sigma_x^2)^{-\frac{m}{2}}.
\]
so altogether (possibly for new constants $C,B$)
\begin{equation}
\label{eqtsattsa2}
\begin{array}{c}
|g_m(S)| 

<C  \sum\limits_{m=0}^{+\infty} B^m \varepsilon^{\frac{m}{2}} (\eta+\frac{\varepsilon\pi}{2}\sigma_x^2)^{-\frac{m}{2}}
 \sum\limits_{|A|=m}   
\frac{1}{\sqrt{A!}} \leqslant \\ { } \\

\leqslant C  \sum\limits_{m=0}^{+\infty} B^m \varepsilon^{\frac{m}{2}} (\eta+\frac{\varepsilon\pi}{2}\sigma_x^2)^{-\frac{m}{2}}
 \sum\limits_{|A|=m}   
\frac{\sqrt{m!}}{\sqrt{A! m!}} \leqslant \\ { } \\

\leqslant C  \sum\limits_{m=0}^{+\infty} B^m \varepsilon^{\frac{m}{2}} (\eta+\frac{\varepsilon\pi}{2}\sigma_x^2)^{-\frac{m}{2}}

\sqrt{ \sum\limits_{|A|=m}   
\binom{m}{A} }

\sqrt{  \frac{  \sum\limits_{|A|=m} 1 }{m!}  }  \leqslant \\ { } \\

\leqslant C  \sum\limits_{m=0}^{+\infty} B^m \varepsilon^{\frac{m}{2}} (\eta+\frac{\varepsilon\pi}{2}\sigma_x^2)^{-\frac{m}{2}}
\sqrt{ n^m  }
\sqrt{
\frac{   n^m }{m!}
} \leqslant C' \sum\limits_{m=0}^{+\infty} \frac{ (B')^m }{ \sqrt{m!} }.
\end{array}
\end{equation}

(For the last step it is important to recall that finally $\eta$ is fixed $\eta=B_0 \varepsilon$ for an appropriate constant $B_0$, see equation (\ref{SetEta})). 

This justifies the interchanging of the integral and series. To proceed to equation (\ref{eqF}) then we carry out the (inverse) Fourier transforms in
\[
\begin{array}{c}
\int{e^{2\pi i Sx-(\eta+\frac{\varepsilon\pi}{2}\sigma_x^2)S^2}\hat{V}(S)
\frac{\varepsilon^m}{m!2^m} [ \sum\limits_{d=1}^n S_d(i\sigma_x^2\partial_{x_d}-\partial_{k_d})]^m
\widetilde{W}_1^\varepsilon(x,k)dS}= \\ { } \\

= \mathcal{F}^{-1}_{S \shortrightarrow x} \left[{ e^{-(\eta+\frac{\varepsilon\pi}{2}\sigma_x^2)S^2}\hat{V}(S)
\frac{\varepsilon^m}{m!2^m} [ \sum\limits_{d=1}^n S_d(i\sigma_x^2\partial_{x_d}-\partial_{k_d})]^m
\widetilde{W}_1^\varepsilon(x,k) }\right]
\end{array}
\]

The proof is complete.

\begin{lemma}[An observation on the Gamma function]
\label{LemGam}
For $m$ large enough
\[
\Gamma \left({ \frac{m-1-\theta}{2} }\right) < (2\pi)^{\frac{1}4}2^{\frac{3+2\theta}{4}} \sqrt{\frac{m!}{2^m m(m-1)}}
\]
\end{lemma}

\noindent {\bf Proof:} The Gamma function satisfies the duplication formula,
\[
\Gamma(z)\Gamma(z+\frac{1}{2})=2^{\frac{1}2-2z} \sqrt{2\pi} \Gamma(2z).
\]
Keeping in mind that $\Gamma(m-1-\theta) < \Gamma(m-1)=(m-2)!=\frac{m!}{m(m-1)}$ (we're working for $m$ large enough), and setting $z=\frac{m-1-\theta}{2}$ it follows that
\[
\begin{array}{c}
\Gamma(\frac{m-1-\theta}{2})\Gamma(\frac{m-\theta}{2})=2^{\frac{1}2-m+1+\theta} \sqrt{2\pi} \Gamma(m-1-\theta) \\ { } \\

\left({\Gamma(\frac{m-1-\theta}{2})}\right)^2<\Gamma(\frac{m-1-\theta}{2})\Gamma(\frac{m-\theta}{2})=2^{\frac{3}2-m+\theta} \sqrt{2\pi} \Gamma(m-1-\theta)<\\ 
\,\,\,\,\,\,\,\,\,\,\,\,\,\,\,\,\,\,\,\,\,\,\,\,\,\,\,\,\,\,\,\,\,\,\,\,\,\,\,\,\,\,\,\,\,\,\,\,\,\,\,\,\,\,\,\,\,\,\,\,\,\,\,\,\,\,\,\,\,\,\,\,\,\,\,\,\,\,\,\,\,\,\,\,
<2^{\frac{3}2+\theta} \sqrt{2\pi} \frac{m!}{2^m m (m-1)}
\end{array}
\]
The proof is complete.

%
%
%

\begin{lemma}[Sobolev norms of delta functions] \label{lemdelsob} Let $z_0 \in \mathbb{R}^{2n}$. Then
\[
2m>n \,\,\, \Rightarrow \,\,\,
\delta(z-z_0) \in H^{-m}(\mathbb{R}^{2n})
\]
\end{lemma}
 
\noindent {\bf Proof:} Take $f \in H^{m}(\mathbb{R}^{2n})$. Then
\[
\begin{array}{c}
| \langle \delta(z-z_0), f(z) \rangle | \leqslant \int{ |\widehat{f}(x,k)|dxdk } =\\ { } \\

= \int{ |\widehat{f}(x,k)| (1+|x|)^{m}(1+|k|)^{m} (1+|x|)^{-m}(1+|k|)^{-m} dk } \leqslant \\ { } \\

\leqslant \sqrt{ \int{ |\widehat{f}(k)|^2 (1+|x|)^{2m}(1+|k|)^{2m} dk } }  \sqrt{ \int{  (1+|k|)^{-2m} dk } \int{  (1+|x|)^{-2m} dk } } \leqslant \\ { } \\

\leqslant C ||f||_{H^m} \left({ 1 + \int\limits_{\rho=1}^{+\infty}{\rho^{n-1-2m} d\rho} }\right)
\end{array}
\]

The proof is complete.

 \begin{appendix}
 \section{Two useful properties of the SWT}\label{app1}
 
 \begin{lemma} \label{thrmUnSmooth} Recall the definition of the spaces $\mathcal{X}^{m,p}$, equation (\ref{eqaloijk}).
Then, for any $m\in \mathbb{R}$, $p\in[1,\infty]$, $s>1$, the following estimate holds
\begin{equation}
|| W-\widetilde{W}||_{\mathcal{X}^{m,p}} \leqslant (1+\frac{\pi}{2}max\{\sigma_x^2,\sigma_k^2 \}) \,\, \varepsilon^{\frac{1}s} \,\, || W ||_{\mathcal{X}^{m+\frac{2}{s},p}}
\end{equation}
\end{lemma}

\noindent {\bf Remark:} The statement and the proof are optimized for $\sigma_x,\sigma_k=O(1)$. When that is not the case, fundamentally the same approach can still be made to work. We don't treat that here explicitly.

\vskip 0.25cm


\noindent {\bf Proof:} For $p\in [1,\infty)$, we have
\begin{equation}
\begin{array}{c}
|| W-\Phi W||_{\mathcal{X}^{m,p}}^p=|| \hat{W}(X,K)(1-e^{-\frac{\varepsilon\pi}{2}(\sigma_x^2X^2+\sigma_k^2K^2)})|(X,K)|^m||_{L^p}^p= \\ { } \\

=\int{ \left|{ (1-e^{-\frac{\varepsilon\pi}{2}[\sigma_x^2X^2+\sigma_k^2K^2]})|(X,K)|^m \hat{W}(X,K) }\right|^p dXdK}=\\ { } \\

=\int\limits_{\sqrt{X^2+K^2}<\frac{1}{\sqrt{\varepsilon}}}{ \left|{ (1-e^{-\frac{\varepsilon\pi}{2}[\sigma_x^2X^2+\sigma_k^2K^2]}) |(X,K)|^m \hat{W}(X,K) }\right|^pdXdK}+\\ 

+\int\limits_{\sqrt{X^2+K^2}>\frac{1}{\sqrt{\varepsilon}}}{ \left|{ (1-e^{-\frac{\varepsilon\pi}{2}[\sigma_x^2X^2+\sigma_k^2K^2]}) |(X,K)|^m \hat{W}(X,K) }\right|^pdXdK} \leqslant \\ { } \\
\end{array}
\end{equation}
\[
\begin{array}{c}
\leqslant \left({ \frac{\pi}{2}max\{\sigma_x^2,\sigma_k^2 \}  }\right)^p \varepsilon^{\frac{p}{s}} \int\limits_{X^2+K^2<\frac{1}{\varepsilon}}{ \left|{ |(X,K)|^{m+\frac{2}{s}} \hat{W}(X,K) }\right|^p dXdK}+\\ 
\,\,\,\,\,\,\,\,\,\,\,\,\,\,\,\,\,\,\,\,\,\,\,\,\,\,\,\,\,\,\,\,\,
+\int\limits_{\sqrt{X^2+K^2}>\frac{1}{\sqrt{\varepsilon}}}{ \left|{|(X,K)|^m  \hat{W}(X,K) }\right|^pdXdK}\leqslant \\ { } \\

\leqslant  \left({ \frac{\pi}{2}max\{\sigma_x^2,\sigma_k^2 \}  }\right)^p \varepsilon^{\frac{p}{s}}
|| W ||_{\mathcal{X}^{m+\frac{2}{s},p}}^p

+\underbrace{\int\limits_{\sqrt{X^2+K^2}>\frac{1}{\sqrt{\varepsilon}}}{ \left|{ |(X,K)|^m \hat{W}(X,K) }\right|^pdXdK}}_{I_2}.
\end{array}
\]
for any $s>1$.

For $I_2$ it suffices to observe that
\begin{equation}
\begin{array}{c}
|| W ||_{\mathcal{X}^{m+\frac{2}{s},p}}^p =  \int{|(X,K)|^{(m+\frac{2}{s})p} |\hat{W}(X,K)|^pdXdK} \geqslant \\ { } \\

\geqslant \varepsilon^{-\frac{p}{s}} \int\limits_{\sqrt{X^2+K^2}>\frac{1}{\sqrt{\varepsilon}}}{|(X,K)|^{mp} |\hat{W}(X,K)|^pdXdK} =\varepsilon^{-\frac{p}{s}} I_2.
\end{array}
\end{equation}

For $p=\infty$ it follows along similar lines (the only difference is technical, we can't raise to the $p=\infty$).

The proof is complete.

It might be instructive to write down explicitly some corollaries of lemma \ref{thrmUnSmooth}:
\begin{eqnarray}
\label{eq5rt67y}
|| W-\widetilde{W}||_{L^2} \leqslant (1+\frac{\pi}{2}max\{\sigma_x^2,\sigma_k^2 \}) \,\, \varepsilon^{\frac{1}2} \,\, || W ||_{H^1}, \\
|| W-\widetilde{W}||_{L^2} \leqslant (1+\frac{\pi}{2}max\{\sigma_x^2,\sigma_k^2 \}) \,\, \varepsilon \,\, || W ||_{H^2}, \\
\label{eq118877hhffdd}
|| W-\widetilde{W}||_{H^{-1}} \leqslant (1+\frac{\pi}{2}max\{\sigma_x^2,\sigma_k^2 \}) \,\, \varepsilon^{\frac{1}2} \,\, || W ||_{L^2}, \\ 
|| W-\widetilde{W}||_{H^{-1}} \leqslant (1+\frac{\pi}{2}max\{\sigma_x^2,\sigma_k^2 \}) \,\, \varepsilon \,\, || W ||_{H^1}.
\end{eqnarray}

\vskip 0.2cm
 
\begin{lemma}[Critical smoothing] \label{lemhuscrit}  If $\sigma_x \sigma_k =1$, then
\[
\widetilde{W}^\varepsilon (x,k)=  \langle  Op_{Weyl} \left({ W^\varepsilon }\right) \phi^\varepsilon_{x,k}, \phi_{x,k}^\varepsilon  \rangle  ,
\]
where
\[
\phi^\varepsilon_{x,k}=\left({ \frac{2}{\varepsilon \sigma_x^2} }\right)^{\frac{n}{4}} e^{-\frac{2\pi i}{\varepsilon} yk - \frac{\pi}{\varepsilon\sigma_x^2}(x-y)^2}.
\]

If $\sigma_x \sigma_k =1$ we will say that we have {\em critical smoothing}, and the SWT coincides with the Husimi transform.

Moreover, if $Op_{Weyl} \left({ W^\varepsilon }\right)$ is a nonnegative definite operator, and $\sigma_x \sigma_k \geqslant 1$, then $\widetilde{W}^\varepsilon \geqslant 0$.
\end{lemma} 
 
\noindent {\bf Proof:} Recall that
\[
\begin{array}{c}
\langle L(x,k), W^\varepsilon[f,g] \rangle= \int { L(x,k) e^{2\pi i k y} \overline{f}(x+\frac{\varepsilon y}{2}) {g}(x-\frac{\varepsilon y}{2})dydxdk }= \\ { } \\

=\varepsilon^{-n} \int{ L(\frac{X+Y}{2},k) e^{2\pi i k \frac{(X-Y)}{\varepsilon}} dk  g(X) dX \overline{f}(Y)dY }=\\ { } \\
= \int{ L(\frac{X+Y}{2},\varepsilon k) e^{2\pi i k (X-Y)} dk  g(X) dX \overline{f}(Y)dY }
 
 =\langle Op_{Weyl}(L) g,f \rangle .
\end{array}
\]
Moreover, it is straightforward to see that, as long as $\sigma_x \sigma_k =1$,
\[
F_{x_0,k_0}(x,k)=\left({ \frac{2}{\varepsilon \sigma_x\sigma_k} }\right)^n {e^{-\frac{2\pi}{\varepsilon} \left[{ \frac{|x-x_0|^2}{\sigma_x^2} + \frac{|k-k_0|^2}{\sigma_k^2} }\right] } }=W^\varepsilon[\phi^\varepsilon_{x_0,k_0}](x,k).
\]
Therefore
\[
\widetilde{W}^\varepsilon(x_0,k_0)= \langle W^\varepsilon , F_{x_0,k_0} \rangle = \langle  Op_{Weyl} \left({ W^\varepsilon }\right) \phi^\varepsilon_{x_0,k_0}, \phi_{x_0,k_0}^\varepsilon  \rangle.
\]

For the second part of the theorem, it is obvious that if $\sigma_x \sigma_k >1$, $\widetilde{W}^\varepsilon$ can be cast as a Gaussian molification of  a Husimi transform. But the nonnegativity of the operator passes on to the Husimi transform, and any Gaussian smoothing preserves it.

 \section{On the initial data of the Wigner equation}\label{app2}

In the case of a pure state $W^\varepsilon[u^\varepsilon](x,k)=\int{ e^{-2\pi i k y}u^\varepsilon(x+\frac{\varepsilon y}{2})\overline{u^\varepsilon}(x-\frac{\varepsilon y}{2})dy }$ the following identity holds:
\[
\int{W^\varepsilon(x,k)dxdk}=||u^\varepsilon||_{L^2}^2.
\]
More generally, in the case of a mixed state $W^\varepsilon(x,k)=\int{ e^{-2\pi i k y}\rho(x+\frac{\varepsilon y}{2},x-\frac{\varepsilon y}{2})dy }$ we can recover the trace of the operator with Weyl symbol $W$ -- equivalently, with kernel $\rho$:
\[
\int{W^\varepsilon(x,k)dxdk}=\int{\rho(x,x)dx}=tr\left({ Op_W W(x,k) }\right).
\]

This means that the conservation of $\int{W^\varepsilon(x,k)dxdk}$ in time evolution reflects a significant physical fact; it emerges therefore as a reasonable normalization to set
\begin{equation}
\label{eqkkaa0jhfc}
\int{W^\varepsilon(x,k)dxdk}=1
\end{equation}
for all $\varepsilon>0$. This is the normalization consistent with the standard scaling of the Schr\"odinger equation, and the concept of a Wigner measure \cite{LP},
\[
W^\varepsilon \rightharpoonup W^0,  \,\,\,\,\,\,\,\,\, \,\,\,\,\,\,\,\,\, \int{W^\varepsilon(x,k)dxdk}=\int{W^0(x,k)dxdk}=1.
\]

\vskip 1cm
At the same time, the $L^2$ norm in also preserved (see theorem \ref{thrmWignL2}), and this is also of physical significance,
\[
|| W^\varepsilon[u^\varepsilon](x,k)||_{L^2}=\varepsilon^{-\frac{n}{2}} ||u^\varepsilon||_{L^2}^2.
\]
In the case of mixed states it corresponds to the Hilber-Schmidt norm,
\[
||W^\varepsilon(x,k)||_{L^2}=\varepsilon^{-\frac{n}{2}} ||\rho(x,y)||_{L^2}= \varepsilon^{-\frac{n}{2}}|| Op_W W(x,k)||_{HS}.
\]

Since we deal with linear homogeneous equations, switching between the two is completely painless. The results here are formulated without any one normalization being fixed, that is why the assumptions involve quantities of the form $\frac{||W_0^\varepsilon||_{H^m}}{||W_0^\varepsilon||_{L^2}}$.

\vskip 1cm
It is instructive to look at pure states and how they scale:

\begin{lemma} \label{lemalan} For the Wigner transform of a pure state, $W^\varepsilon=W^\varepsilon[u^\varepsilon]$, the following holds: 
\[
\begin{array}{c}
||\partial_{x_i} W^\varepsilon||_{L^2} \leqslant 2 \varepsilon^{-\frac{n}{2}} \,\, ||\partial_{x_i} u^\varepsilon||_{L^2}\,\, || u^\varepsilon||_{L^2}, \\ { } \\

||\partial_{k_i} W^\varepsilon||_{L^2} \leqslant 4\pi \varepsilon^{-\frac{n}{2}} \,\, ||\frac{x_i}{\varepsilon} u^\varepsilon||_{L^2}\,\, || u^\varepsilon||_{L^2}.
\end{array}
\]
\end{lemma}

\noindent {\bf Proof: } Denote by $Q$ the mapping
\[
Q: f(x,y) \mapsto \int{e^{2\pi i ky} f(x+\frac{\varepsilon y}{2},x-\frac{\varepsilon y}{2})  dy}.
\]
For example, the Wigner transform of a pure state can be recast as
\[
W^\varepsilon[u^\varepsilon]=Q[ u^\varepsilon(x)\overline{u}^\varepsilon(y) ].
\]
It is easy to check that
\begin{equation}
\label{eqwwwnnnn22}
|| Q f ||_{L^2} = \varepsilon^{-\frac{n}{2}} ||f||_{L^2}.
\end{equation}
Moreover, it is easy to check that
\[
\begin{array}{c}
\partial_{x_i} W^\varepsilon[u^\varepsilon] = Q\left[{(\partial_{x_i}+\partial_{y_i}) u^\varepsilon(x)\overline{u}^\varepsilon(y) }\right], \\ { } \\

\partial_{k_i} W^\varepsilon[u^\varepsilon] = Q\left[{-2\pi i\frac{x_i-y_i}{\varepsilon} u^\varepsilon(x)\overline{u}^\varepsilon(y) }\right].
\end{array}
\]
The proof is complete.

\vskip 0.75 cm
The importance of the following lemma is that it assures us that  for coherent states the convergence to the Wigner measure is in fact in $H^{-r}$ norm, not just in weak-$*$ sense. This is in particular useful when working with strong semiclassical limits for pure states.

\begin{lemma}[Wigner measure for coherent states] \label{LEMATARA} Assume $u^\varepsilon_0(x)$ is a coherent state,
\begin{equation}
u^\varepsilon_0(x)=\varepsilon^{-\frac{n}{4}} a(\frac{x-x_0}{\sqrt{\varepsilon}}) e^{2\pi i k_0 x},
\end{equation}
with $a \in \mathcal{S}$.

Then, for $\phi \in \mathcal{S}$,  and for all $r \in \mathbb{N}$, $2r>n$, there is a constant $C$ depending on $a(x)$, $r$ such that
\[
|\langle W^\varepsilon[u^\varepsilon_0] - \delta(x_0,k_0),\phi \rangle | \leqslant  C || \phi||_{H^{r+2}} \left({ \varepsilon^{\frac{2r-n}{4n}} + \varepsilon^{\frac{1}{2n}} }\right).
\]
\end{lemma}

\noindent {\bf Remark: } Similar results can be shown for WKB initial data, i.e. $u^\varepsilon(x)=A(x)e^{\frac{2\pi i}{\varepsilon}S(x)}$. 

\vskip 0.2cm
\noindent {\bf Proof: } It suffices to work for $x_0=k_0=0$. Indeed,
\[
\begin{array}{c}
\langle W^\varepsilon[u^\varepsilon_0],\phi \rangle = \varepsilon^{-\frac{n}{2}}\int{ e^{-2\pi i k y}
 a(\frac{x+\frac{\varepsilon y}{2}}{\sqrt{\varepsilon}}) 
a(\frac{x-\frac{\varepsilon y}{2}}{\sqrt{\varepsilon}}) \phi(x,k)dydxdk}= \\ { } \\

=  \int{ e^{-2\pi i k y} a(x+\frac{\sqrt{\varepsilon} y}{2})  a(x-\frac{\sqrt{\varepsilon} y}{2}) \phi(\sqrt{\varepsilon}x,k)dydxdk}= \\ { } \\

=  \underbrace{\int{ e^{-2\pi i k y} a(x+\frac{\sqrt{\varepsilon} y}{2})  a(x-\frac{\sqrt{\varepsilon} y}{2}) \phi(0,k)dydxdk} }_{I} + \\ { } \\
+\underbrace{ \int{ e^{-2\pi i k y} a(x+\frac{\sqrt{\varepsilon} y}{2})  a(x-\frac{\sqrt{\varepsilon} y}{2}) R_1 dydxdk} }_{I_1} ,
\end{array}
\]
where of course $R_1$ is the remainder of the Taylor expansion. We will collect and estimate all the errors in the end; in the meantime we have
\[
\begin{array}{c}
\int{ e^{-2\pi i k y} a(x+\frac{\sqrt{\varepsilon} y}{2})  a(x-\frac{\sqrt{\varepsilon} y}{2}) \phi(0,k)dydxdk}= \\ { } \\
=\int{ e^{-2\pi i k y} \phi(0,k) dk a(x+\frac{\sqrt{\varepsilon} y}{2})  a(x-\frac{\sqrt{\varepsilon} y}{2}) dydx}= \\ { } \\

=\int{  \widehat{\phi}_2(0,y) a(x+\frac{\sqrt{\varepsilon} y}{2})  a(x-\frac{\sqrt{\varepsilon} y}{2}) dydx}=
\int{  \widehat{\phi}_2(0,y) [a(x)+E_1] [a(x)+E_1']  dydx}=\\ { } \\

= \phi(0,0)+ \int{\widehat{\phi}_2(0,y) [a(x)E_1+a(x)E_1'+E_1E_1']dydx},
\end{array}
\]
where of course $E_1$, $E_1'$ are the Taylor remainders from Taylor expanding $a(x\pm \frac{\sqrt{\varepsilon}y}{2})$ around $x$, and $\widehat{\phi}_2(x,y)=\int{e^{-2\pi i k y}\phi(x,k)dk}$.

Now, about the errors: $|R_1| \leqslant |x| \sqrt{\varepsilon}  || \nabla \phi ||_{L^\infty}$, and therefore
\[
\begin{array}{c}
|I_1| \leqslant \sqrt{\varepsilon}  || \nabla \phi ||_{L^\infty} |\int{ e^{-2\pi i k y} a(x+\frac{\sqrt{\varepsilon} y}{2})  a(x-\frac{\sqrt{\varepsilon} y}{2}) |x| dydxdk} | \leqslant \\ { } \\

\leqslant \sqrt{\varepsilon}  || \nabla \phi ||_{L^\infty} |\int{ \delta(y-0) a(x+\frac{\sqrt{\varepsilon} y}{2})  a(x-\frac{\sqrt{\varepsilon} y}{2}) |x| dydx} | \leqslant 

\sqrt{\varepsilon}  || \nabla \phi ||_{L^\infty} \int{a^2(x) |x|dx} .
\end{array}
\]
Moreover, $|E_1|,|E_1'| \leqslant \frac{1}{2} |y|\sqrt{\varepsilon} ||\nabla a ||_{L^\infty}$. Therefore,
\[
\begin{array}{c}
|\int{\widehat{\phi}_2(0,y) a(x)E_1dydx}| \leqslant \frac{\sqrt{\varepsilon} ||\nabla a ||_{L^\infty}}{2}|\int{\widehat{\phi}_2(0,y) a(x)  |y| dydx}|= \\ { } \\

= \frac{\sqrt{\varepsilon} ||\nabla a ||_{L^\infty}}{2}\int{|\widehat{\phi}_2(0,y)|\,\,|y| dy} \int{ |a(x)|  dx}.
\end{array}
\]

Finally, for the $E_1E_1'$ term, one observes that
\[
\begin{array}{c}
\int{\widehat{\phi}_2(0,y) E_1E_1'dydx} = \\ { } \\
=\int{\widehat{\phi}_2(0,y) [a(x+\frac{\sqrt{\varepsilon} y}{2})  a(x-\frac{\sqrt{\varepsilon} y}{2})-a(x)E_1-a(x)E_1'-a^2(x)]dydx},
\end{array}
\]
and since we can select a $M>1$ such that
\[
\int\limits_{|x|>M, \, y \in \mathbb{R}^n}{\widehat{\phi}_2(0,y) [a(x)E_1+a(x)E_1'+a^2(x)]dydx} \leqslant \sqrt{\varepsilon} \int\limits_{ y \in \mathbb{R}^n}{\widehat{\phi}_2(0,y)dy},
\]
it follows that
\[
\begin{array}{c}
|\int{\widehat{\phi}_2(0,y) E_1E_1'dydx}| \leqslant |\int\limits_{|x|<M, \, y \in \mathbb{R}^n}{\widehat{\phi}_2(0,y) E_1E_1'dydx}| +\\ { } \\

+ \sqrt{\varepsilon} \int\limits_{ y \in \mathbb{R}^n}{|\widehat{\phi}_2(0,y)|dy} + | \int\limits_{|x|>M, \, y \in \mathbb{R}^n}{\widehat{\phi}_2(0,y) a(x+\frac{\sqrt{\varepsilon} y}{2})  a(x-\frac{\sqrt{\varepsilon} y}{2})dydx} |.
\end{array}
\]
Now we have
\[
\begin{array}{c}
|\int\limits_{|x|<M, \, y \in \mathbb{R}^n}{\widehat{\phi}_2(0,y) E_1E_1'dydx}| \leqslant \frac{\varepsilon}{4}  ||\nabla a ||^2_{L^\infty} \int\limits_{|x|<M} {dx} \int{\widehat{\phi}_2(0,y) |y|^2dy},
\end{array}
\]
\[
\begin{array}{c}
I_2:=| \int\limits_{|x|>M, \, y \in \mathbb{R}^n}{\widehat{\phi}_2(0,y) a(x+\frac{\sqrt{\varepsilon} y}{2})  a(x-\frac{\sqrt{\varepsilon} y}{2})dydx} |\leqslant \\ { } \\

\leqslant | \int\limits_{|x|>M, \, |y|>M }{\widehat{\phi}_2(0,y) a(x+\frac{\sqrt{\varepsilon} y}{2})  a(x-\frac{\sqrt{\varepsilon} y}{2})dydx} |+\\

+| \int\limits_{|x|>M, \, |y|<M }{\widehat{\phi}_2(0,y) a(x+\frac{\sqrt{\varepsilon} y}{2})  a(x-\frac{\sqrt{\varepsilon} y}{2})dydx} | \leqslant \\ { } \\

\leqslant |\int\limits_{ |y|>M }{\widehat{\phi}_2(0,y) \left({\int{ a(x+\frac{\sqrt{\varepsilon} y}{2})  a(x-\frac{\sqrt{\varepsilon} y}{2})dx} }\right)dy} |+ \\

+| \int\limits_{|x|>M, \, |y|<M }{\widehat{\phi}_2(0,y) a(x+\frac{\sqrt{\varepsilon} y}{2})  a(x-\frac{\sqrt{\varepsilon} y}{2})dydx} |.
\end{array}
\]
At this point observe that\footnote{ by restricting, without loss of generality, $\varepsilon \in (0,1)$}
\[
|x|>M, |y|<M \,\, \Rightarrow \,\, \frac{1}{2}|x| \leqslant  |x \pm \frac{\sqrt{\varepsilon}}{2}y| \leqslant 2|x|;
\]
now making use of the fact that $a \in \mathcal{S}$ it follows that, as long as $|x|>M, |y|<M$, there are $C_2>0, s_1 > n$ such that
\[
|a(x+\frac{\sqrt{\varepsilon} y}{2})  a(x-\frac{\sqrt{\varepsilon} y}{2})| \leqslant C_2 \frac{1}{(1+|x|)^{s_1}}
\]
so that finally
\[
\begin{array}{c}
I_2 \leqslant ||a(x)||_{L^2}^2 \int\limits_{ |y|>M }{|\widehat{\phi}_2(0,y)|  dy} 

+C_2  \int\limits_{|y|<M }{|\widehat{\phi}_2(0,y)|dy} \int\limits_{|x|>M} {\frac{dx}{(1+|x|)^{s_1}}}
\end{array}
\]

To summarize (and keeping all the $M$ dependencies explicit -- that will be important later),
\begin{equation}
\label{eqapql8xne}
\begin{array}{c}
\left|{ \langle W^\varepsilon[u^\varepsilon_0]-\delta(x-0,k-0),\phi \rangle }\right| \leqslant \\ { } \\

\leqslant  ||a(x)||_{L^2}^2 \int\limits_{ |y|>M }{|\widehat{\phi}_2(0,y)|  dy} 

+C_2  \int\limits_{|y|<M }{|\widehat{\phi}_2(0,y)|dy}\,\, M^{n-s_1} + \\ { } \\

+C \varepsilon ||\nabla a ||^2_{L^\infty} M^n \int{\widehat{\phi}_2(0,y) |y|^2dy}

+\frac{\sqrt{\varepsilon} ||\nabla a ||_{L^\infty}}{2}\int{|\widehat{\phi}_2(0,y)|\,\,|y| dy} \int{ |a(x)|  dx}.
\end{array}
\end{equation}
Obviously
\begin{equation}
\label{eqapql8xne1}
\begin{array}{c}
\int\limits_{|y|<M }{|\widehat{\phi}_2(0,y)|dy} \leqslant \int{|\widehat{\phi}(x,y)|dxdy}, \\ { } \\

\int{\widehat{\phi}_2(0,y) |y|^2dy} \leqslant \int{|\widehat{\phi}(x,y)|\,\,|y|^2 dxdy}, \\ { } \\

\int{|\widehat{\phi}_2(0,y)|\,\,|y| dy} \leqslant \int{|\widehat{\phi}(x,y)|\,\,|y| dxdy}.
\end{array}
\end{equation}
Now it suffices to see that
\begin{equation}
\label{eqapql8xne2}
\begin{array}{c}
\int{|\widehat{\phi}(x,y)| (|1+y|)^2 dxdy} =\\ { } \\

= \int{|\widehat{\phi}(x,y)| (|1+y|)^2 \frac{1}{(|1+x|)^r(|1+y|)^r}(|1+x|)^r(|1+y|)^r dxdy} \leqslant \\ { } \\

\leqslant C || \phi ||_{H^{r+2}} || \frac{1}{(|1+x|)^r(|1+y|)^r} ||_{L^2} = C' || \phi ||_{H^{r+2}}
\end{array}
\end{equation}
for any $2r>n$.

To complete the estimation, one observes that
\begin{equation}
\begin{array}{c}
\int\limits_{ |y|>M }{|\widehat{\phi}_2(0,y)|  dy} \leqslant \int\limits_{|y|>M} {|\widehat{\phi}(x,y)|dxdy}= \\ { } \\

=\int\limits_{|y|>M} {|\widehat{\phi}(x,y)| \frac{1}{(|1+x|)^r(|1+y|)^r}(|1+x|)^r(|1+y|)^r dxdy} \leqslant \\ { } \\

\leqslant C ||\phi||_{L^2} \sqrt{ \int\limits_{|y|>M}{ \frac{1}{(|1+x|)^{2r}(|1+y|)^{2r}} dxdy  }} \leqslant \\ { } \\

\leqslant C' ||\phi||_{L^2} \sqrt{ M^{n-2r} }.
\end{array}
\end{equation}

So by setting $M=\varepsilon^{-\frac{1}{2n}}$, $s_1=n+1$ in equation (\ref{eqapql8xne}), and applying the estimations (\ref{eqapql8xne1}), (\ref{eqapql8xne2}), the result follows.

\end{appendix}

\section*{Acknowledgment}

A. G. Athanassoulis would like to thank the  D\'epartement de Math\'ematiques et Applications of the \'Ecole Normale Superieure for financial support.


\end{document}